\documentclass{amsart}
\usepackage{amscd,amssymb,amsxtra}

\theoremstyle{plain}
\newtheorem{Thm}{Theorem}[section]
\newtheorem{Prop}[Thm]{Proposition}
\newtheorem{Lem}[Thm]{Lemma}
\newtheorem{Cor}[Thm]{Corollary}
\newtheorem{Claim}[Thm]{Claim}

\theoremstyle{definition}
\newtheorem{Def}[Thm]{Definition}
\newtheorem{Question}[Thm]{Question}

\theoremstyle{remark}

\numberwithin{equation}{section}

\newcommand{\dom}{\operatorname{dom}}
\newcommand{\res}{\upharpoonright}
\newcommand{\diag}{\operatorname{diag}}

\begin{document}

\title{Infinite products of finite simple groups}
\author{Jan Saxl}
\address[Jan Saxl]%
{DPMMS \\
16 Mill Lane \\
Cambridge CB2 1SB \\
England}
\author{Saharon Shelah}
\address[Saharon Shelah]%
{Mathematics Department \\
The Hebrew University \\
Jerusalem \\
Israel}
\thanks{The research of the second author was partially supported
by the BSF. Publication 584 of the second author.}
\author{Simon Thomas}
\address[Simon Thomas]%
{Mathematics Department \\
Bilkent University \\
Ankara \\
Turkey}
\address[Saharon Shelah and Simon Thomas]%
{Mathematics Department \\
Rutgers University \\
New Brunswick \\
New Jersey 08903}
\thanks{The research of the third author was partially supported
by NSF Grants.}

\begin{abstract}
We classify the sequences $\langle S_{n} \mid n \in \mathbb{N} \rangle$
of finite simple nonabelian groups such that $\prod_{n} S_{n}$ has
uncountable cofinality.
\end{abstract}
\maketitle
\section{Introduction} \label{S:intro}
Suppose that $G$ is a group that is not finitely generated. Then $G$
can be expressed as the union of a chain of proper subgroups. The
cofinality of $G$, written $c(G)$, is defined to be the least cardinal
$\lambda$ such that $G$ can be expressed as the union of a chain of
$\lambda$ proper subgroups. Groups of uncountable cofinality were
first considered by Serre in his study of groups acting on trees.

\begin{Def}{\cite[p.58]{se}} \label{D:FA}
A group $H$ has property $(FA)$ if and only if whenever $H$ acts
without inversion on a tree $T$, then there exists a vertex
$t \in T$ such that $h(t) = t$ for all $h \in H$.
\end{Def}

In \cite{se}, Serre characterised the groups which have 
property (FA).

\begin{Thm}{\cite{se}} \label{T:FA}
The group $H$ has property (FA) if and only if the following
three conditions are satisfied.
\begin{enumerate}
\item[(1)] $H$ is not a nontrivial free product with amalgamation.
\item[(2)] $\mathbb{Z}$ is not a homomorphic image of $H$.
\item[(3)] If $H$ is not finitely generated, then $c(H) > \omega$.
\end{enumerate}
\end{Thm}

This result led to the question of whether there exist any
natural examples of uncountable groups with property (FA).
Let $\langle G_{n} \mid n \in \mathbb{N} \rangle$ be a sequence
of nontrivial finite groups. Then $\prod_{n} G_{n}$ denotes
the full direct product of the groups $G_{n}$, $n \in \mathbb{N}$.
By Bass \cite{ba}, if $H$ is a profinite group and $H$ acts
without inversion on the tree $T$, then for every $h \in H$
there exists $t \in T$ such that $h(t) = t$. This implies
that $H$ satisfies conditions \ref{T:FA}(1) and \ref{T:FA}(2).
In particular, we see that the profinite group
$\prod_{n} G_{n}$ has property (FA) if and only if
$c(\prod_{n} G_{n}) > \omega$. The following result, which
was proved by Koppelberg and Tits, provided the first examples
of uncountable groups with property (FA).

\begin{Thm}{\cite{kt}} \label{T:tits}
Let $F$ be a nontrivial finite group and let $G_{n} = F$ for
all $n \in \mathbb{N}$. Then $c(\prod_{n} G_{n}) > \omega$ if and
only if $F$ is perfect.
\end{Thm}

Suppose that $F$ is perfect.
Since $\left| \prod_{n} G_{n} \right| = 2^{\omega}$, Theorem 
\ref{T:tits} yields that
\[
\omega_{1} \leq c \left(\prod_{n} G_{n} \right) \leq 2^{\omega}.
\]
This suggests the problem of trying to compute the exact value
of $c(\prod_{n} G_{n})$. (Of course, this problem is only
interesting if $2^{\omega} > \omega_{1}$.) The following
result is an immediate consequence of a theorem of Koppelberg
\cite{ko}.

\begin{Thm} \label{T:exact}
If $F$ is a nontrivial finite perfect group and $G_{n} = F$ for all
$n \in \mathbb{N}$, then $c(\prod_{n} G_{n}) = \omega_{1}$.
\end{Thm}

\begin{proof}
If $\langle g(n) \rangle_{n} \in \prod_{n} G_{n}$ and $\pi \in F$,
let $X_{\pi}(g) = \{ n \in \mathbb{N} \mid g(n) = \pi \}$. Then
$\{ X_{\pi}(g) \mid \pi \in F \}$ yields a partition of $\mathbb{N}$
into finitely many pieces. Consider the powerset $\mathcal{P}(\mathbb{N})$
as a Boolean algebra. By Koppelberg \cite{ko}, we can express
\[
\mathcal{P}(\mathbb{N}) = \underset{\alpha < \omega_{1}}{\bigcup} B_{\alpha}
\]
as the union of a chain of $\omega_{1}$ proper Boolean subalgebras.
For each $\alpha < \omega_{1}$, define
\[
H_{\alpha} = \{ g \in \prod_{n} G_{n} \mid
X_{\pi}(g) \in B_{\alpha} \text{ for all } \pi \in F \}.
\]
Then it is easily checked that $H_{\alpha}$ is a proper subgroup
of $\prod_{n} G_{n}$. Clearly $\prod_{n} G_{n} =
\underset{\alpha < \omega_{1}}{\bigcup} H_{\alpha}$, and so
$c(\prod_{n} G_{n}) \leq \omega_{1}$.
\end{proof}

The above results suggest the following questions.

\begin{Question} \label{Q:cof}
For which sequences $\langle S_{n} \mid n \in \mathbb{N} \rangle$
of finite simple nonabelian groups, do we have that
$c(\prod_{n} S_{n}) > \omega$?
\end{Question}

\begin{Question} \label{Q:exact}
Suppose that $\langle S_{n} \mid n \in \mathbb{N} \rangle$ is a 
sequence of finite simple nonabelian groups such that
$c(\prod_{n} S_{n}) > \omega$. Is it possible to compute
the exact value of $c(\prod_{n} S_{n})$?
\end{Question}

It may be helpful to give a word of explanation concerning
Question \ref{Q:exact}. The point is that it may be impossible
to compute the exact value of $c(\prod_{n} S_{n})$ in $ZFC$.
For example, consider the group $Sym(\mathbb{N})$ of all
permutations of $\mathbb{N}$. In \cite{mn}, Macpherson and Neumann
showed that $c(Sym(\mathbb{N})) > \omega$. Later Sharp and Thomas
\cite{st1} proved that it is consistent that $c(Sym(\mathbb{N}))$
and $2^{\omega}$ can be any two prescribed regular uncountable
cardinals subject only to the requirement that
$c(Sym(\mathbb{N})) \leq 2^{\omega}$. Hence it is impossible to
compute the exact value of $c(Sym(\mathbb{N}))$ in $ZFC$. (The 
theorem of Macpherson and Neumann suggests that $Sym(\mathbb{N})$
is probably another natural example of an uncountable group
with property (FA). In the final section of this paper, we
shall confirm that this is true.)

The following result shows that there exist sequences
$\langle S_{n} \mid n \in \mathbb{N} \rangle$ of finite simple nonabelian 
groups such that $c(\prod_{n} S_{n}) = \omega$.

\begin{Thm} \label{T:count}
Let $\langle S_{n} \mid n \in \mathbb{N} \rangle$ be a sequence of
finite simple nonabelian groups. Suppose that there exists an
infinite subset $I$ of $\mathbb{N}$ such that the following conditions
are satisfied.
\begin{enumerate}
\item[(1)] There exists a fixed (possibly twisted) Lie type $L$ such
that for all $n \in I$, $S_{n} = L(q_{n})$ for some prime power
$q_{n}$.
\item[(2)] If $n$, $m \in I$ and $n < m$, then $q_{n} < q_{m}$.
\end{enumerate}
Then $c(\prod_{n} S_{n}) = \omega$.
\end{Thm}

Here $L(q_{n})$ denotes the group of Lie type $L$ over the finite
field $GF(q_{n})$. The proof of Theorem \ref{T:count} makes use of
the following easy observation.

\begin{Lem} \label{L:normal}
Suppose that $N \vartriangleleft G$ and that $G/N$ is not finitely
generated. Then $c(G) \leq c(G/N)$.
\end{Lem}
\begin{flushright}
$\square$
\end{flushright}

\begin{proof}[Proof of Theorem \ref{T:count}]
By Lemma \ref{L:normal}, we can suppose that $I = \mathbb{N}$. Let
$\mathcal{D}$ be a nonprincipal ultrafilter on $\mathbb{N}$, and
let $N$ be the set of elements $g = \langle g(n) \rangle_{n}
\in \prod_{n} S_{n}$ such that
$\{ n \in \mathbb{N} \mid g(n) = 1 \} \in \mathcal{D}$. Then $N$ is
a normal subgroup of $\prod_{n} S_{n}$, and $\prod_{n} S_{n} /N$
is the ultraproduct $G = \prod_{n} S_{n} / \mathcal{D}$. (See Section
9.5 of Hodges \cite{h}.) By Lemma \ref{L:normal}, it is enough to
show that $c(G) = \omega$.

There exists a fixed integer $d$ such that each of the groups $L(q_{n})$
has a faithful $d$-dimensional linear representation over the field
$GF(q_{n})$. Since the class of groups with a faithful $d$-dimensional
linear representation is first-order axiomatisable, it follows that $G$
has a faithful $d$-dimensional linear representation over some field
$K$. (For example, see Section 6.6 of Hodges \cite{h}. It is perhaps
worth mentioning that every known proof only yields the existence of
a set of axioms for this class. The problem of finding an explicit
intelligible set of axioms remains open.)
To simplify
notation, we shall suppose that $G \leqslant GL(d,K)$. We also suppose that
$K$ has been chosen so that $G \cap GL(d,K^{\prime})$ is a proper subgroup
of $G$ for every proper subfield $K^{\prime}$ of $K$. By Exercise 9.5.5 of
Hodges \cite{h}, $|G| = 2^{\omega}$. It follows that $|K| = 2^{\omega}$,
and hence $K$ has transcendence dimension $2^{\omega}$ over its
prime subfield $k$. Let $B$ be a transcendence basis of $K$ over $k$.
Express $B = \underset{n < \omega}{\bigcup} B_{n}$ as the union
of a chain of proper subsets. For each $n < \omega$, let $K_{n}$
be the algebraic closure of $B_{n}$ in $K$. Then each $K_{n}$ is
a proper subfield of $K$, and 
$K = \underset{n < \omega}{\bigcup} K_{n}$. For each $n < \omega$,
let $G_{n} = G \cap GL(d,K_{n})$. Then $G = 
\underset{n < \omega}{\bigcup} G_{n}$, and each $G_{n}$ is a proper
subgroup of $G$. Hence $c(G) = \omega$.
\end{proof}

The main result of this paper is that the converse of Theorem
\ref{T:count} is also true.

\begin{Thm} \label{T:class}
Suppose that $\langle S_{n} \mid n \in \mathbb{N} \rangle$ is a 
sequence of finite simple nonabelian groups such that
$c(\prod_{n} S_{n}) = \omega$. Then there exists an infinite
subset $I$ of $\mathbb{N}$ such that conditions \ref{T:count}(1) and 
\ref{T:count}(2) are satisfied.
\end{Thm}

Now suppose that $\langle S_{n} \mid n \in \mathbb{N} \rangle$ is a sequence of 
finite simple nonabelian groups such that $c(\prod_{n} S_{n}) > \omega$.
If there exists an infinite subset $J$ of $\mathbb{N}$ such that 
$S_{n} = S_{m}$ for all $n$, $m \in J$, then Lemma \ref{L:normal}
and Theorem \ref{T:exact} imply that $c(\prod_{n} S_{n}) = \omega_{1}$.
This is the only case in which we have been able to compute the exact
value of $c(\prod_{n} S_{n})$ in $ZFC$.

\begin{Question} \label{Q:unc}
Is it consistent that there exists a sequence
$\langle S_{n} \mid n \in \mathbb{N} \rangle$ of finite simple nonabelian
groups such that $c(\prod_{n} S_{n}) > \omega_{1}$?
\end{Question}

We hope that Question \ref{Q:unc} has a positive answer, as this
would lead to some very attractive problems. For example, consider
the following question. (We suspect that it cannot be answered in 
$ZFC$.)

\begin{Question} \label{Q:dream}
Is $c(\prod_{n} Alt(n+5)) = c( \prod_{n} PSL(n+3,2))$?
\end{Question}

In Section \ref{S:con}, we shall prove the following consistency result.
Amongst other things, it shows that it is impossible to prove in $ZFC$
that $c(Sym(\mathbb{N})) = c(\prod_{n} S_{n})$ for some sequence
$\langle S_{n} \mid n \in \mathbb{N} \rangle$ of finite simple nonabelian
groups.

\begin{Thm} \label{T:con}
It is consistent that both of the following statements are true.
\begin{enumerate}
\item[(1)] $c(Sym(\mathbb{N})) = \omega_{2} = 2^{\omega}$.
\item[(2)] $c(\prod_{n} G_{n}) \leq \omega_{1}$ for {\em every\/}
sequence $\langle G_{n} \mid n \in \mathbb{N} \rangle$ of nontrivial
finite groups.
\end{enumerate}
\end{Thm}

The following problem is also open. (Of course, a negative answer
to Question \ref{Q:unc} would yield a negative answer to
Question \ref{Q:con}.)

\begin{Question} \label{Q:con}
Is it consistent that there exists a sequence
$\langle S_{n} \mid n \in \mathbb{N} \rangle$ of finite simple nonabelian
groups such that $c(\prod_{n} S_{n}) > c(Sym(\mathbb{N}))$?
\end{Question}

This paper is organised as follows. In Section \ref{S:alt}, we
shall prove that if $\langle S_{n} \mid n \in \mathbb{N} \rangle$
is a sequence of finite alternating groups, then
$c(\prod_{n} S_{n}) > \omega$. In Section \ref{S:lin}, we shall
prove Theorem \ref{T:class} in the special case when each $S_{n}$
is a projective special linear group. In Section \ref{S:class},
we shall complete the proof of Theorem \ref{T:class}. Our proof
makes use of the fact that there are only finitely many sporadic
finite simple groups, and thus relies on the classification of
the finite simple groups. Section
\ref{S:con} contains the proof of Theorem \ref{T:con}. In
Section \ref{S:sym}, we shall prove that $Sym(\mathbb{N})$ has
property (FA).

Our notation is standard, but a couple of points should be mentioned.
Suppose that $G$ is a subgroup of $Sym(\Omega)$. If each nonidentity
element $g \in G$ is fixed-point-free, then $G$ is said to act
{\em semiregularly\/} on $\Omega$. If $G$ acts transitively and
semiregularly, then $G$ is said to act {\em regularly\/} on
$\Omega$. In this paper, permutation groups and linear groups
always act on the left. Thus, for example, we have that
\[
(\, 1 \, 2 \, 3 \,)(\, 1 \, 3 \, 5 \, 7 \,)(\, 1 \, 2 \, 3 \,)^{-1}
=(\, 2 \, 1 \, 5 \, 7\,).
\]
We follow the usual convention of regarding each ordinal as the set
of its predecessors. Thus $\omega = \mathbb{N}$. Also if $a$, $b$ are
natural numbers such that $a > b$, then their set-theoretic
difference is $a \smallsetminus b =
\{ b, b+1, \dots , a-1 \}$. If $A$ is a matrix, then $A^{T}$ denotes
the transpose of $A$.

\section{Infinite products of alternating groups} \label{S:alt}
In this section, we shall prove the following special case of 
Theorem \ref{T:class}.

\begin{Thm} \label{T:alt}
Let $\langle S_{n} \mid n \in \mathbb{N} \rangle$ be a sequence of finite
simple nonabelian groups. If each $S_{n}$ is an alternating group,
then $c(\prod_{n} S_{n}) > \omega$.
\end{Thm}

We shall make use of the following two results, which will be
used repeatedly throughout this paper.

\begin{Prop}{\cite{th}} \label{P:gen}
Suppose that $G$ is not finitely generated and that $H$ is a subgroup
of $G$. If $G$ is finitely generated over $H$, then $c(H) \leq c(G)$.
\end{Prop}

\begin{proof}
Let $c(G) = \lambda$. Express $G = \underset{\alpha < \lambda}{\bigcup}
G_{\alpha}$
as the union of a chain of $\lambda$ proper subgroups. Let
$H_{\alpha} = H \cap G_{\alpha}$. Then $H = \underset{\alpha < \lambda}
{\bigcup} 
H_{\alpha}$. Since $G$ is finitely generated over $H$, each $H_{\alpha}$
is a proper subgroup of $H$. Thus $c(H) \leq \lambda$.
\end{proof}

\begin{Prop} \label{P:finite}
Let $\langle S_{n} \mid n \in \mathbb{N} \rangle$ be a sequence of nontrivial
finite perfect groups. Suppose that there exists a finite set
$\mathcal{F}$ of groups such that $S_{n} \in \mathcal{F}$ for all
$n \in \mathbb{N}$. Then $c(\prod_{n} S_{n}) > \omega$.
\end{Prop}

\begin{proof}
By Proposition \ref{P:gen}, we can suppose that for each $S \in \mathcal{F}$,
the set \\
$\{ n \in \mathbb{N} \mid S_{n} = S \}$ is either infinite or
empty. Since the class of groups of uncountable cofinality is
closed under taking finite direct products, Theorem \ref{T:exact} implies
that $c(\prod_{n} S_{n}) > \omega$.
\end{proof}

We shall begin the proof of Theorem \ref{T:alt}
by making a couple of easy reductions. For each $m \geq 5$, let
$P_{m} = \prod_{n} S^{m}_{n}$, where $S^{m}_{n} = Alt(m)$ for
all $n \in \mathbb{N}$. Let $G_{0} = \prod_{m \geq 5}P_{m}$. Then
Lemma \ref{L:normal} implies that it is enough to prove that
$c(G_{0}) > \omega$. Let $G_{1} = \prod_{m \geq 8}P_{m}$. Then
$G_{0} = P_{5} \times P_{6} \times P_{7} \times G_{1}$. By Theorem \ref{T:exact},
$c(P_{m}) = \omega_{1}$ for all $m \geq 5$. Hence it is enough to prove
that $c(G_{1}) > \omega$. Finally let $G_{2} = \prod_{m \geq 3}P_{2^{m}}$.
Then Theorem \ref{T:alt} is an immediate consequence of the following
two results.

\begin{Lem} \label{L:red}
$c(G_{1}) = c(G_{2})$.
\end{Lem}

\begin{Thm} \label{T:main}
$c(G_{2}) > \omega$.
\end{Thm}

First we shall prove Lemma \ref{L:red}. Note that Lemma \ref{L:normal}
implies that $c(G_{1}) \leq c(G_{2})$. Our proof that $c(G_{2}) \leq
c(G_{1})$ is based upon Proposition \ref{P:gen}.

Let $I = \{ \langle m,n \rangle \mid 8 \leq m \in \mathbb{N}, n \in \mathbb{N} \}$.
Then $G_{1} = \underset{ \langle m,n \rangle \in I}{\prod}S^{m}_{n}$,
where $S^{m}_{n} = Alt(m)$. For each 
$\langle m,n \rangle \in I$, let $t$ be the integer such that
$2^{t} \leq m < 2^{t+1}$ and let
$T^{m}_{n} = Alt(2^{t}) \leqslant S^{m}_{n}$. Then we can identify
$G_{2}$ with the subgroup
$\underset{\langle m,n \rangle \in I}{\prod} T^{m}_{n}$ of $G_{1}$.
By Proposition \ref{P:gen}, it is enough to prove the following
result.

\begin{Lem} \label{L:gen}
$G_{1}$ is finitely generated over $G_{2}$.
\end{Lem}

This is the first of the many places in this paper where we need to 
prove that an infinite product of groups is finitely generated over
an infinite product of subgroups. A moment's thought shows that
such results require ``uniform generation'' results for the 
corresponding sequences of groups. We shall make repeated use of
the following easy observation.

\begin{Prop} \label{P:word}
Let $\langle H_{n} \mid n \in \mathbb{N} \rangle$ and 
$\langle G_{n} \mid n \in \mathbb{N} \rangle$ be sequences of groups
such that $H_{n} \leqslant G_{n}$ for all $n \in \mathbb{N}$. Suppose
that there exists a word $w(x_{1}, \dots , x_{s}, y_{1}, \dots , y_{t})$
from the free group on $\{ x_{1}, \dots , x_{s}, y_{1}, \dots , y_{t} \}$
such that the following condition is satisfied.
\begin{enumerate}
\item[(\ref{P:word})] For all $n \in \mathbb{N}$, there exist
elements $\theta_{1}, \dots , \theta_{t} \in G_{n}$ such that each
$\phi \in G_{n}$ can be expressed as
$\phi = w(\psi_{1}, \dots , \psi_{s}, \theta_{1},
\dots , \theta_{t})$
for some $\psi_{1}, \dots , \psi_{s} \in H_{n}$.
\end{enumerate}
Then there exist elements $g_{1}, \dots , g_{t} \in \prod_{n} G_{n}$
such that $\prod_{n} G_{n} =
\langle \prod_{n} H_{n}, g_{1}, \dots , g_{t} \rangle$.
\end{Prop}
\begin{flushright}
$\square$
\end{flushright}

Lemma \ref{L:gen} is a consequence of the
following ``uniform generation'' results for the finite alternating
groups, which will also be needed in the proof of Theorem \ref{T:main}.

\begin{Lem} \label{L:uni1}
Let $m \geq 3$ and let $\theta = (\ m-2\ m-1\ )(\ m\ m+1\ ) \in Alt(m+1)$.
Then for every $\phi \in Alt(m+1)$, there exists
$\psi_{1}$, $\psi_{2}$, $\psi_{3} \in Alt(m)$ such that
\[
\phi = \psi_{1} \theta \psi_{2} \theta \psi_{3}.
\]
\end{Lem}

\begin{proof}
If $\phi \in Alt(m)$, then we can take $\psi_{1} = \phi$ and
$\psi_{2} = \psi_{3} = id$. So suppose that 
$\phi \in Alt(m+1) \smallsetminus Alt(m)$. Let
$\psi_{2} = (\ m-2\ m-1\ m\ )$. Then
$\tau = \theta \psi_{2} \theta = (\ m-1\ m-2\ m+1\ )$. Since $Alt(m+1)$
acts 2-transitively on $\{ 1, \dots , m+1 \}$, we have the double
coset decomposition
\[
Alt(m+1) = Alt(m) \cup Alt(m) \tau Alt(m).
\]
Thus $\phi \in Alt(m) \tau Alt(m)$, and so there exist
$\psi_{1}$, $\psi_{3} \in Alt(m)$ such that
$\phi = \psi_{1} \theta \psi_{2} \theta \psi_{3}$.
\end{proof}

Let $m = 4n$ for some $n \geq 2$. 
A permutation $\pi \in Alt(m)$ is said to have type $2^{2n}$
if $\pi$ is the product of $2n$ disjoint transpositions. Thus
$\pi^{2} = id$ and $\pi$ is fixed point free. The set of permutations
$\pi \in Alt(m)$ of type $2^{2n}$ forms a single conjugacy class in
$Alt(m)$.

\begin{Lem}[Brenner \cite{br}] \label{L:brenner}
Let $m = 4n$ for some $n \geq 2$. Let $C$ be the conjugacy class
of $Alt(m)$ consisting of all permutations of type $2^{2n}$. Then
for every $\phi \in Alt(m)$, there exist
$\pi_{1}, \dots , \pi_{4} \in C$ such that
$\phi = \pi_{1} \dots \pi_{4}$.
\end{Lem}
\begin{flushright}
$\square$
\end{flushright}

\begin{Lem} \label{L:uni2}
Suppose that $m = 8n$ for some $n \geq 1$. Let $\Delta_{0} = \{1, \dots , 4n \}$,
$\Delta_{1} = \{4n+1, \dots ,8n \}$ and let $\Gamma =
Alt(\Delta_{0}) \times Alt(\Delta_{1})$. Let
$\theta = \prod_{i=1}^{4n}(\ i\ 2i\ )(\ 4n+i\ 2i-1\ )$. 
Then every
$\phi \in Alt(m)$ can be expressed as a product
\[
\phi = \psi_{1} \theta \psi_{2} \theta \psi_{3} \theta \psi_{4}
\theta \psi_{5} \theta \psi_{6} \theta \psi_{7} \theta \psi_{8}
\theta \psi_{9}
\]
for some $\psi_{1}, \dots , \psi_{9} \in
\Gamma$.
\end{Lem}

\begin{proof} 
By Lemma \ref{L:brenner}, it is
enough to show that each
permutation $\phi \in Alt(m)$ of type $2^{4n}$ can be expressed as a product
$\phi = \psi_{1} \theta \psi_{2} \theta \psi_{3}$
for some $\psi_{1}$, $\psi_{2}$ , $\psi_{3} \in \Gamma$.
Let $A = \{\ell \in \Delta_{0} \mid \phi( \ell) \in \Delta_{0} \}$
and $B = \{\ell \in \Delta_{1} \mid \phi( \ell) \in \Delta_{1} \}$.
Then $|A| = |B| = 2s$ for some $0 \leq s \leq 4n$. Let
$C = \{1, \dots ,8n \} \smallsetminus (A \cup B)$. Then
$|C \cap \Delta_{0}| = |C \cap \Delta_{1}| = 4n-2s$. Let $t = 4n-2s$.
Let $\Delta_{2} = \{2i \mid 1 \leq i \leq 4n \}$ and 
$\Delta_{3} = \{2i-1 \mid 1 \leq i \leq 4n \}$.

\noindent
{\em Case 1.} Suppose that $t \geq 2n$. Choose a subset
$D \subseteq C \cap \Delta_{0}$ of size $2n$, and let
$E = \phi[D]$. Then there exists $\psi_{1} \in \Gamma$ such that
$\psi_{1}[D] = \Delta_{2} \cap \Delta_{0}$ and
$\psi_{1}[E] = \Delta_{2} \cap \Delta_{1}$. This implies that
\[
\psi_{1} \phi \psi_{1}^{-1} \in Alt(\Delta_{2}) \times Alt(\Delta_{3})
= \theta \left( Alt(\Delta_{0}) \times Alt(\Delta_{1}) \right) \theta.
\]
Thus we have that
\[
\psi_{2} = \theta \psi_{1} \phi \psi_{1}^{-1} \theta
\in \Gamma
\]
and so
\[
\phi = \psi_{1}^{-1} \theta \psi_{2} \theta \psi_{1}
\]
is a suitable product.

\noindent
{\em Case 2.} Suppose that $t < 2n$. Then $s > n$. Choose
$\phi$-invariant subsets $D \subseteq \Delta_{0}$ and
$E \subseteq \Delta_{1}$ such that $|D| = |E| = 2n$. Then
there exists $\psi_{1} \in \Gamma$ such that
$\psi_{1}[D] = \Delta_{2} \cap \Delta_{0}$ and
$\psi_{1}[E] = \Delta_{2} \cap \Delta_{1}$. Arguing as in
Case 1, we see that there exists $\psi_{2} \in \Gamma$ such
that $\phi = \psi_{1}^{-1} \theta \psi_{2} \theta \psi_{1}$. 
\end{proof}

\begin{proof}[Proof of Lemma \ref{L:gen}]
For each $0 \leq i \leq 7$, let
\[
H_{i} = \prod\{S_{n}^{m} \mid \langle m,n \rangle \in I,
m \equiv i \text{ (mod 8) } \}
\]
and
\[
K_{i} = \prod\{T_{n}^{m} \mid \langle m,n \rangle \in I,
m \equiv i \text{ (mod 8) } \}.
\]
Then $G_{1} = \prod_{i=0}^{7} H_{i}$ and
$G_{2} = \prod_{i=0}^{7} K_{i}$. Clearly it is enough to show
that $H_{i}$ is finitely generated over $K_{i}$ for each
$0 \leq i \leq 7$.

First consider the case when $i=0$. Let 
$\langle m,n \rangle \in I$ satisfy $m = 8s$ for some $s \geq 1$,
and let $t$ be the integer such that $2^{t} \leq m < 2^{t+1}$.
Define $\Delta_{0}^{m,n} = \{1, \dots ,4s \}$ and
$\Delta_{1}^{m,n} = \{4s+1, \dots ,8s \}$. Then
$Alt(\Delta_{0}^{m,n}) \leqslant T_{n}^{m} = Alt(2^{t})$.
There exists an element $\phi \in Alt(m) = S_{n}^{m}$ such
that $\phi Alt(\Delta_{0}^{m,n}) \phi^{-1}
= Alt(\Delta_{1}^{m,n})$. Hence there exists $g_{1} \in H_{0}$
such that
\[
\prod \{ Alt(\Delta_{0}^{m,n}) \times Alt(\Delta_{1}^{m,n}) \mid
\langle m,n \rangle \in I, m \equiv 0 \text{ (mod 8) } \}
\leqslant \langle K_{0}, g_{1} \rangle .
\]
Now Lemma \ref{L:uni2} implies that there exists an element
$g_{2} \in H_{0}$ such that 
$H_{0} = \langle K_{0}, g_{1}, g_{2} \rangle$.

Next consider the case when $i=1$. For each 
$\langle m,n \rangle \in I$ with $m \equiv 1$ (mod 8), let
$U_{n}^{m} =Alt(m-1) \leqslant S_{n}^{m}$. By the previous
paragraph, there exist elements $g_{1}$, $g_{2} \in H_{1}$
such that
\[
\prod \{ U_{n}^{m} \mid \langle m,n \rangle \in I,
m \equiv 1 \text{ (mod 8) } \} \leqslant 
\langle K_{1}, g_{1}, g_{2} \rangle .
\]
Now Lemma \ref{L:uni1} implies that there exists an element
$g_{3} \in H_{1}$ such that 
$H_{1} = \langle K_{1}, g_{1}, g_{2}, g_{3} \rangle$. Continuing
in this fashion, we can successively deal with the remaining
cases. 
\end{proof}

The rest of this section is devoted to the proof of Theorem \ref{T:main}.
Suppose that $c(G_{2}) = \omega$. Express
$G_{2} = \underset{t < \omega}{\bigcup}H_{t}$ as the union of a chain of
$\omega$ proper subgroups. Our strategy will be to define by induction
on $t < \omega$
\begin{enumerate}
\item a sequence of elements $f_{t} \in G_{2}$;
\item a strictly increasing sequence of integers $i_{t}$
such that $f_{t} \in H_{i_{t}}$;
\item a sequence of elements $g_{t} \in G_{2} \smallsetminus H_{i_{t}}$.
\end{enumerate}
These sequences will be chosen so that there exists an element
$h \in G_{2}$ such that $hf_{t}h^{-1} = g_{t}$ for all $t < \omega$.
But this implies that $h \notin \underset{t < \omega}{\bigcup}H_{t}$,
which is the desired contradiction.

Let $J = \{ \langle m,n \rangle \mid 3 \leq m \in \mathbb{N}, n \in \mathbb{N} \}$.
Then $G_{2} = \underset{\langle m,n \rangle \in J}{\prod}A_{n}^{m}$,
where $A_{n}^{m} = Alt(2^{m})$ for all $n \in \mathbb{N}$. The elements
$g_{t} = \langle g_{t}(m,n) \rangle_{m,n} \in \prod_{m,n} A_{n}^{m}$,
$t < \omega$, will be chosen so that for each 
$\langle m,n \rangle \in J$, the sequence
\[
g_{0}(m,n),\ g_{1}(m,n), \dots ,\ g_{t}(m,n)
\]
is a generic sequence of elements of $Alt(2^{m})$, in the following
sense.

\begin{Def} \label{D:generic}
If $0 \leq t \leq m-1$, then the sequence
$\pi_{0}, \dots , \pi_{t}$ of elements of $Alt(2^{m})$
is a {\em generic sequence\/} if
\begin{enumerate}
\item the subgroup $\langle \pi_{0}, \dots ,\pi_{t} \rangle$ is
elementary abelian of order $2^{t+1}$;
\item if $id \ne \phi \in \langle \pi_{0}, \dots ,\pi_{t} \rangle$,
then $\phi$ is a permutation of type $2^{2^{m-1}}$. (In other words,
$\langle \pi_{0}, \dots ,\pi_{t} \rangle$ acts semiregularly on
$\{1, \dots ,2^{m} \}$.)
\end{enumerate}
If $m-1 \leq t < \omega$, then the sequence $\pi_{0}, \dots ,\pi_{t}$
of elements of $Alt(2^{m})$ is a generic sequence if
\begin{enumerate}
\item[(a)] $\pi_{0}, \dots ,\pi_{m-1}$ is a generic sequence;
\item[(b)] $\pi_{\ell} = \pi_{m-1}$ for all $m-1 \leq \ell \leq t$.
\end{enumerate}
\end{Def}

It is an easy exercise to show that for each $t < \omega$, there
exists a unique generic sequence $\pi_{0}, \dots ,\pi_{t}$ in
$Alt(2^{m})$ up to conjugacy within $Sym(2^{m})$; and two such
generic sequences up to conjugacy within $Alt(2^{m})$ if
$m \geq 3$ and $2 \leq t < \omega$. (We shall not make
use of this observation in the proof of Theorem \ref{T:main}.)

To begin the induction, choose any element 
$f_{0} = \langle f_{0}(m,n) \rangle_{m,n} \in
G_{2} = \prod_{m,n} A_{n}^{m}$ such that $f_{0}(m,n)$ is a 
permutation of type $2^{2^{m-1}}$ for each 
$\langle m,n \rangle \in J$; and let $i_{0}$ be an integer
such that $f_{0} \in H_{i_{0}}$. Let $f_{0}^{G_{2}}$ be the
conjugacy class of $f_{0}$ in $G_{2}$. Then Lemma \ref{L:brenner}
implies that $G_{2} = \langle f_{0}^{G_{2}} \rangle$. Hence
there exists an element $h_{0} \in G_{2}$ such that
$g_{0} = h_{0} f_{0} h_{0}^{-1} \notin H_{i_{0}}$.
Now suppose that $t \geq 0$
and that we have defined 
\begin{enumerate}
\item a sequence of elements $f_{j} \in G_{2}$,
\item a strictly increasing sequence of integers $i_{j}$ 
such that $f_{j} \in H_{i_{j}}$,
\item a sequence of elements $g_{j} \in G_{2} \smallsetminus H_{i_{j}}$, and
\item a sequence of elements $h_{j} \in G_{2}$
\end{enumerate}
for $0 \leq j \leq t$ such that the following conditions hold.
\begin{enumerate}
\item[(a)] $f_{0}(m,n), \dots , f_{t}(m,n)$ is a generic sequence in
$A_{n}^{m} = Alt(2^{m})$ for all $m,n$.
\item[(b)] If $0 \leq j \leq  k \leq t$, then
$h_{k} f_{j} h_{k}^{-1} = g_{j}$.
\item[(c)] If $m-1 \leq j \leq t$ and $n \in \mathbb{N}$, then
$h_{j}(m,n) = h_{m-1}(m,n)$.
\end{enumerate}  

First we shall define $f_{t+1}$.

\noindent
{\em Case 1.} Suppose that $m-1 \leq t$ and $n \in \mathbb{N}$.
Then we define $f_{t+1}(m,n) = f_{m-1}(m,n)$.

\noindent
{\em Case 2.} Suppose that $t < m-1$ and $n \in \mathbb{N}$.
We shall set up some notation which will be used during the rest of
this section. Let
\[
E(m,n) = \langle g_{0}(m,n), \dots ,g_{t}(m,n) \rangle .
\]
Then $E(m,n)$ is an elementary abelian group of order $2^{t+1}$
acting semiregularly on $\{1, \dots ,2^{m} \}$. Let
\[
\{1, \dots ,2^{m} \} = \Phi^{m,n}_{1} \cup \dots \cup \Phi^{m,n}_{2^{m-t-1}}
\]
be the decomposition into $E(m,n)$-orbits. Then $E(m,n)$ acts regularly
on $\Phi^{m,n}_{i}$ for each $1 \leq i \leq 2^{m-t-1}$. Choose 
$\alpha_{i} \in \Phi^{m,n}_{i}$ for each $1 \leq i \leq 2^{m-t-1}$.
Let $E(m,n) = \{ \pi_{k} \mid 1 \leq k \leq 2^{t+1} \}$,
where $\pi_{1} = id$, and define
\[
\Delta^{m,n}_{k} = \{ \pi_{k}(\alpha_{i}) \mid 1 \leq i \leq 2^{m-t-1} \}
\]
for each $1 \leq k \leq 2^{t+1}$. Then the diagonal subgroup
\begin{align*}
D(m,n) &= Diag(Alt(\Delta^{m,n}_{1}) \times \dots \times Alt(\Delta^{m,n}_{2^{t+1}}))
\\
       &= \{ \prod_{i=1}^{2^{t+1}} \pi_{i} \phi \pi_{i} \mid
\phi \in Alt(\Delta^{m,n}_{1}) \}
\end{align*}
is contained in the centraliser of $E(m,n)$ in $Alt(2^{m})$. Let
$\tau (m,n) \in D(m,n)$ be any permutation of type $2^{2^{m-1}}$, 
and define
\[
f_{t+1}(m,n) = h_{t}(m,n)^{-1} \tau (m,n) h_{t}(m,n).
\]
This completes the definition of
$f_{t+1} = \langle f_{t+1}(m,n) \rangle_{m,n}$.

Next we choose $i_{t+1}$ to be an integer such that
\begin{enumerate}
\item[(i)] $i_{t} < i_{t+1}$ and $f_{t+1} \in H_{i_{t+1}}$; and
\item[(ii)] $P^{*} = \prod \{ A^{m}_{n} \mid 3 \leq m \leq t+3, 
n \in \mathbb{N} \} \leqslant H_{i_{t+1}}$.
\end{enumerate}
(Proposition \ref{P:finite} implies that
$c(P^{*}) > \omega$. Hence $i_{t+1}$ can be chosen so that (ii) also
holds.)

Finally we shall define $h_{t+1}$ and $g_{t+1}$.

\noindent
{\em Case 1.} Suppose that $m-1 \leq t$ and $n \in \mathbb{N}$. Then we
define $h_{t+1}(m,n) = h_{m-1}(m,n)$ and $g_{t+1}(m,n) = g_{m-1}(m,n)$.

\noindent
{\em Case 2.} Suppose that $t < m-1$ and $n \in \mathbb{N}$. Then we choose
a suitable element $\sigma(m,n) \in D(m,n)$ and define
\[
h_{t+1}(m,n) = \sigma(m,n) h_{t}(m,n)
\]
and
\begin{align*}
g_{t+1}(m,n) &= h_{t+1}(m,n) f_{t+1}(m,n) h_{t+1}(m,n)^{-1} \\
             &= \sigma(m,n) \tau(m,n) \sigma(m,n)^{-1}.
\end{align*}
Of course, a suitable choice means one such that
$g_{t+1} = \langle g_{t+1}(m,n) \rangle_{m,n} \notin H_{i_{t+1}}$. This
completes the successor step of the induction, {\em provided\/} that
a suitable choice exists.

\begin{Claim} \label{C:ind}
There exists a choice of $\sigma (m,n)$ for $m > t+1$
and $n \in \mathbb{N}$  such that $g_{t+1} \notin H_{i_{t+1}}$.
\end{Claim}

\begin{proof}
Suppose that for every choice of the sequence
\[
\langle \sigma (m,n) \mid t< m-1, n \in \mathbb{N} \rangle
\]
we have that $g_{t+1} \in H_{i_{t+1}}$. Then we shall prove that
$G_{2}$ is finitely generated over $H_{i_{t+1}}$. But this means that there
exists $r \in \mathbb{N}$ with $i_{t+1} \leq r$ such that $H_{r} = G_{2}$, which
is a contradiction.

Let $J^{\prime} = \{ \langle m,n \rangle \mid t+4 \leq m \in \mathbb{N},
n \in \mathbb{N} \}$ and let
$P^{\prime} = \prod \{A^{m}_{n} \mid \langle m,n \rangle \in J^{\prime} \}$.
Thus $G_{2} = P^{*} \times P^{\prime}$. Note that if 
$\pi(m,n) \in D(m,n)$ is any element of type $2^{2^{m-1}}$, then there
exists $\sigma (m,n) \in D(m,n)$ such that
$\sigma (m,n) \tau (m,n) \sigma (m,n)^{-1} = \pi (m,n)$. Using the 
fact that $P^{*} \leqslant H_{i_{t+1}}$, we see that the following statement
holds.
\begin{enumerate}
\item[($\dag$)] Suppose that $\pi = \langle \pi (m,n) \rangle_{m,n} \in P^{\prime}$.
If $\pi (m,n) \in D(m,n)$ is an element of type $2^{2^{m-1}}$ for all
$\langle m,n \rangle \in J^{\prime}$, then $\pi \in H_{i_{t+1}}$.
\end{enumerate}
Using Lemma \ref{L:brenner}, we see that
$\prod \{ D(m,n) \mid \langle m,n \rangle \in J^{\prime} \}
\leqslant H_{i_{t+1}}$.
Now let $\theta_{1} = \langle \theta_{1} (m,n) \rangle_{m,n}
\in P^{\prime}$ be an element such that
$\theta_{1} (m,n) \in Alt(\Delta^{m,n}_{1})$ is a permutation of type
$2^{2^{m-t-2}}$ for each $\langle m,n \rangle \in J^{\prime}$. Then
$\{ \psi \theta_{1} (m,n) \psi^{-1} \mid \psi \in D(m,n) \}$ is the
conjugacy class in $Alt(\Delta^{m,n}_{1})$ of all permutations of
type $2^{2^{m-t-2}}$. Using Lemma \ref{L:brenner} again, we see
that
\[
\prod \{ Alt(\Delta^{m,n}_{1}) \mid \langle m,n \rangle \in J^{\prime} \}
\leqslant \langle H_{i_{t+1}}, \theta_{1} \rangle.
\]
Continuing in this fashion, we find that there
exist $\theta_{1}, \dots ,\theta_{2^{t+1}} \in P^{\prime}$
such that
\[
\prod \{ Alt(\Delta^{m,n}_{1}) \times \dots \times Alt(\Delta^{m,n}_{2^{t+1}})
\mid \langle m,n \rangle \in J^{\prime} \} \leqslant
\langle H_{i_{t+1}}, \theta_{1}, \dots ,\theta_{2^{t+1}} \rangle.
\]
By repeatedly applying Lemma \ref{L:uni2}, we now see that there exists
a finite subset $F$ of $P^{\prime}$ such that
\[
P^{\prime} \leqslant \langle H_{i_{t+1}}, \theta_{1}, \dots ,
\theta_{2^{t+1}}, F \rangle.
\]
Hence $G_{2} = P^{*} \times P^{\prime}$ is finitely generated over
$H_{i_{t+1}}$, which is a contradiction.
\end{proof}

Thus the induction can be carried out for all $t < \omega$. Define
the element
$h = \langle h(m,n) \rangle_{m,n} \in G_{2}$ by
$h(m,n) = h_{m-1}(m,n)$. Then we have that 
$h f_{t} h^{-1} = g_{t}$ for all $t < \omega$, which is a contradiction.
This completes the proof of Theorem \ref{T:main}.

\section{Infinite products of special linear groups} \label{S:lin}
In this section, we shall prove the following result.

\begin{Thm} \label{T:lin}
Suppose that $\langle SL(d_{n},q_{n}) \mid n \in \mathbb{N} \rangle$ is a sequence
of finite special linear groups which satisfies the following conditions.
\begin{enumerate}
\item[(1)] If $d_{n} = 2$, then $q_{n} > 3$.
\item[(2)] There does {\em not\/} exist an infinite subset $I$ of 
$\mathbb{N}$ and an integer $d$ such that
\begin{enumerate}
\item $d_{n} = d$ for all $n \in I$; and
\item if $n$, $m \in I$ and $n < m$, then $q_{n} < q_{m}$.
\end{enumerate}
\end{enumerate}
Then $c(\prod_{n} SL(d_{n},q_{n})) > \omega$.
\end{Thm}

Using Theorem \ref{T:lin} and Lemma \ref{L:normal}, we see that
Theorem \ref{T:class} is true in the special case when each
$S_{n}$ is a projective special linear group.

Our strategy in the proof of Theorem \ref{T:lin} will be the same
as that in the proof of Theorem \ref{T:alt}. We shall begin by
defining the notion of a generic sequence of elements in
$SL(2^{m},q)$. Let $V = V(n,q)$ be an n-dimensional vector space
over $GF(q)$, and let $\mathcal{B}$ be a basis of $V$. Then $Sym(\mathcal{B})$
denotes the group of permutation matrices with respect to the 
basis $\mathcal{B}$. Note that for any finite field $GF(q)$, we have
that $Alt(\mathcal{B}) \leqslant SL(n,q)$.

\begin{Def} \label{D:lgen}
If $0 \leq t \leq m-1$, then the sequence $\pi_{0}, \dots , \pi_{t}$
of elements of $SL(2^{m},q)$ is a {\em generic sequence\/} if there
exists a basis $\mathcal{B}$ of $V(2^{m},q)$ such that
\begin{enumerate}
\item the group $\langle \pi_{0}, \dots , \pi_{t} \rangle$ is an
elementary abelian subgroup of $Alt(\mathcal{B})$ of order $2^{t+1}$;
\item $\langle \pi_{0}, \dots , \pi_{t} \rangle$ acts semiregularly
on $\mathcal{B}$.
\end{enumerate}
If $m-1 \leq t < \omega$, then the sequence $\pi_{0}, \dots , \pi_{t}$
of elements of $SL(2^{m},q)$ is a generic sequence if
\begin{enumerate}
\item[(a)] $\pi_{0}, \dots , \pi_{m-1}$ is a generic sequence;
\item[(b)] $\pi_{\ell} = \pi_{m-1}$ for all $m-1 \leq \ell \leq t$.
\end{enumerate}
\end{Def}

First we shall prove an analogue of Lemma \ref{L:brenner}. Let
$m = 4n$ for some $n \geq 1$. Then $C(m,q)$ denotes the conjugacy
class in $SL(m,q)$ consisting of all elements $\pi$ such that
$\pi$ is represented by 
a permutation matrix of type $2^{2n}$ with respect to
some basis $\mathcal{B}$ of $V(m,q)$. (It is easily checked that the
set of such elements forms a single conjugacy class in $SL(m,q)$.)
Note that $SL(m,q)$ has a maximal torus $T_{1}$ of order
$(q^{4n}-1)/(q-1)$. By Zsigmondy's theorem \cite{zs}, there exists a primitive
prime divisor $p > 2$ of $q^{4n}-1$. Let $\psi \in T_{1}$ be an element
of order $p$. 
Then $\psi$ is clearly a regular element of $T_{1}$. (A semisimple element
$g \in SL(m,q)$ is {\em regular\/} if and only if it lies in a unique
maximal torus.)

\begin{Lem} \label{L:reg}
With the above hypotheses, there exist elements $\pi_{1}$,
$\pi_{2} \in C(m,q)$ such that $\psi = \pi_{1} \pi_{2}$.
\end{Lem}

\begin{proof}
Regard $K = GF(q^{4n})$ as a $4n$-dimensional vector space over
$GF(q)$. Let $\tau \in K$ generate a normal basis of $K$ over
$GF(q)$; ie. $\mathcal{B} =
\{ \tau, \tau^{q}, \tau^{q^{2}}, \dots , \tau^{q^{4n-1}} \}$
is a basis of $K$ over $GF(q)$. Let $f \in Aut(K)$ be the
Frobenius automorphism; so that $f(\alpha) = \alpha^{q}$ for
all $\alpha \in K$. Let $g = f^{2n}$. Then $g$ is represented
by a permutation matrix of type $2^{2n}$ with respect to the
basis $\mathcal{B}$. Thus $g \in C(m,q)$.

Let $K^{*} = \langle \alpha \rangle$, and let 
$\beta = \alpha^{q^{2n}-1}$. Then $\beta$ has order $q^{2n}+1$,
and $g(\beta) = \beta^{-1}$. Consider the primitive prime
divisor $p$ of $q^{4n}-1 = (q^{2n}-1)(q^{2n}+1)$. Then clearly
$p$ divides $q^{2n}+1$. Thus there exists an element
$\gamma \in \langle \beta \rangle$ of order $p$; and we can
suppose that $\psi \in T_{1}$ is the linear transformation
defined by $\psi(x) = \gamma x$ for all $x \in K$.
Since $g(\gamma) = \gamma^{-1}$, we see that
$g \psi g^{-1} = \psi^{-1}$. Since $\psi$ has odd order, the
involutions $g$ and $g \psi$ are conjugate in the dihedral
group $\langle g, g \psi \rangle$. The result follows.
\end{proof}

\begin{Thm} \label{T:saxl}
Suppose that $m = 4n$ for some $n \geq 1$. Then for every 
$\phi \in SL(m,q)$, there exist $\pi_{1}, \dots ,\pi_{10}
\in C(m,q)$ such that $\phi = \pi_{1} \dots \pi_{10}$.
\end{Thm}

\begin{proof}
Let $G = SL(4n,q)$ and let $\psi \in T_{1}$ be as above. Let
$\tau \in G$ be an element of order $q^{4n-1}-1$ and let $T_{2}$
be the maximal torus which contains $\tau$. (Of course, $\tau$
is a regular element of $T_{2}$.) 
Let $C_{1}$, $C_{2}$ be the 
conjugacy classes of $\psi$, $\tau$ repectively. We claim that
the product $C_{1}C_{2}$ of these two classes covers all of
$G \smallsetminus Z(G)$. Using Lemma \ref{L:reg}, this implies
that each element of $C_{2}$ is a product of 3 elements of
$C(m,q)$; and hence every element of $G \smallsetminus Z(G)$
is a product of 5 elements of $C(m,q)$. The result follows.

The proof of the claim uses character theory and follows
\cite[pp. 96--99]{msw} very closely. For any conjugacy class
$C_{3}$ of $G$ and $\sigma \in C_{3}$, define
\[
m(C_{1},C_{2},C_{3}) =
\frac{ \left| G \right|^{2} }{ \left| C_{G}(\psi) \right|
\left| C_{G}(\tau) \right| \left| C_{G}(\sigma) \right| }
\sum \frac{ \chi(\psi) \chi(\tau) \chi(\sigma) }
{\chi(1)}
\]
where the summation runs over the irreducible characters $\chi$
of $G$. By a well-known class formula, $m(C_{1},C_{2},C_{3})$ is
equal to the number of triples $(a_{1},a_{2},a_{3})$ such that
$a_{i} \in C_{i}$ and $a_{1}a_{2}a_{3} = 1$. It therefore suffices
to show that the character sum involved in the formula for
$m(C_{1},C_{2},C_{3})$ is positive for any class $C_{3}$ of
non-central elements of $G$.

Now the values of the irreducible characters of $G$ on semisimple
elements can be calculated from the values of the Deligne-Lusztig
characters. (See \cite[Chapter 7]{ca2}.) The Deligne-Lusztig characters
$R_{T,\theta}$ are parametrized by pairs $(T,\theta)$. The equivalence
relation of geometric conjugacy on these pairs yields a partition of
the irreducible characters of $G$ into disjoint series as follows.
The geometric conjugacy classes of pairs $(T,\theta)$ can be
parametrized by the conjugacy classes $(s)$ of semisimple elements
in the dual group $\hat{G} = PGL(4n,q)$. Let $\mathcal{E}(s)$ be the
set of irreducible characters occurring as a constituent in one of 
the $R_{T,\theta}$ with $(T,\theta)$ corresponding to $(s)$. Then
the sets $\mathcal{E}(s)$, where $(s)$ runs over the set of conjugacy
classes of semisimple elements of $\hat{G}$, form a partition of
the set of irreducible characters of $G$. ( See
\cite[7.3.8 and 7.5.8]{ca2}.)

The $R_{T,\theta}$ span the space of class functions restricted to
semisimple elements. (See \cite[7.5.7]{ca2}.) In particular, suppose
that $\rho$ is a semisimple element of $G$ and that $\chi \in 
\mathcal{E}(s)$. Then if $R_{T,\theta}(\rho) = 0$ for all pairs 
$(T,\theta)$ corresponding to $(s)$, we have that $\chi(\rho) = 0$.
Now $R_{T,\theta} = R(s)$ vanishes on the regular elements of the 
torus $T^{\prime}$ if the element $s$ is not conjugate in $\hat{G}$
to an element of the dual $\hat{T^{\prime}}$ of $T^{\prime}$. Hence
the $R_{T,\theta}$ not vanishing on either of the classes $C_{1}$,
$C_{2}$ will correspond to semisimple classes $(s)$ in $\hat{G}$
such that $s \in \hat{T_{1}} \cap \hat{T_{2}}$. Let 
$\bar{T_{i}}$ be the preimage of $\hat{T_{i}}$ in $GL(4n,q)$.
Then $\bar{T_{1}}$ is cyclic of order $q^{4n}-1$, and $\bar{T_{2}}$
is the product of two cyclic groups of orders $q^{4n-1}-1$ and $q-1$.
Furthermore, both $\bar{T_{1}}$ and $\bar{T_{2}}$ contain
$Z(GL(4n,q))$. Since $\left( \left| \bar{T_{1}} \right|, 
\left| \bar{T_{2}} \right|
\right) = (q-1)(4n,q-1)$, it follows that $\hat{T_{1}} \cap \hat{T_{2}} 
=1$. Thus we need only consider the set $\mathcal{E}(1)$ of unipotent
characters of $G$. It is well-known that if the degree of an 
irreducible character is divisible by the full power of a prime $r$
dividing the order of a group, then the character vanishes on all
$r$-singular elements of the group. Using this result, an inspection
of the degrees of the unipotent characters of $G$ in 
\cite[p. 465]{ca2} shows that only two irreducible characters
contribute to the character sum in the above formula; namely, the
principal character and the Steinberg character $St$. It follows
that
\[
m(C_{1},C_{2},C_{3}) =
\frac{ \left| G \right|^{2} }{ \left| C_{G}(\psi) \right|
\left| C_{G}(\tau) \right| \left| C_{G}(\sigma) \right| }
\left( 1 + \frac{ St(\psi) St(\tau) St(\sigma)}{St(1)} \right)
\]
and so \cite[6.4]{ca2}
\[
m(C_{1},C_{2},C_{3}) =
\frac{ \left| G \right|^{2} }{ \left| C_{G}(\psi) \right|
\left| C_{G}(\tau) \right| \left| C_{G}(\sigma) \right| }
\left( 1 - \frac{ St(\sigma) }{ St(1)} \right).
\]
This is 0 precisely when $\sigma \in Z(G)$, as claimed. 
\end{proof}

Next we shall prove the analogue of Lemma \ref{L:uni1}. It is
easier to state the result in terms of infinite products of
groups, rather than in terms of ``uniform generation''. Let
$\langle SL(d_{n}, q_{n}) \mid n \in \mathbb{N} \rangle$ 
be a sequence of special linear groups.
Fix some $n \in \mathbb{N}$. Let $SL(d_{n},q_{n})$ act on the vector space
$V(d_{n},q_{n})$ in the natural manner. Extend this action to
$V(d_{n}+1,q_{n}) = V(d_{n},q_{n}) \oplus
\langle v_{d_{n}+1} \rangle$ by specifying that
$\pi(v_{d_{n}+1}) = v_{d_{n}+1}$ 
for all $\pi \in SL(d_{n},q_{n})$. Using this
extended action, we can regard $SL(d_{n},q_{n})$ as a subgroup
of $SL(d_{n}+1,q_{n})$.

\begin{Lem} \label{L:step}
$\prod_{n} SL(d_{n}+1,q_{n})$ is finitely generated over
$\prod_{n} SL(d_{n},q_{n})$.
\end{Lem}

\begin{proof}
We shall make use of the Bruhat decomposition
\[
SL(d,q) = \underset{w \in W}{\bigcup} BwB
\]
of the special linear group, where $B$ is a Borel subgroup and
$W$ is the Weyl group. Fix some integer $n \in \mathbb{N}$. Choose a basis
$\{ v_{i} \mid 1 \leq i \leq d_{n} \}$ of $V(d_{n},q_{n})$. We 
shall regard each element of $SL(d_{n}+1,q_{n})$ as a matrix with
respect to the basis $\mathcal{B} = \{ v_{i} \mid 1 \leq i \leq d_{n}+1 \}$.
Note that we have identified $SL(d_{n},q_{n})$ with the subgroup
\[
S_{n} = \left \{
\begin{pmatrix}
A & \boldsymbol{0} \\
\boldsymbol{0} & 1
\end{pmatrix}
\mid A \in SL(d_{n},q_{n}) \right \}
\]
of $SL(d_{n}+1,q_{n})$. Define
\[
T_{n} = \left \{
\begin{pmatrix}
1 & \boldsymbol{0} \\
\boldsymbol{0} & A 
\end{pmatrix}
\mid A \in SL(d_{n},q_{n}) \right \} .
\]
Then there exists $\pi \in SL(d_{n}+1,q_{n})$ such that
$\pi S_{n} \pi^{-1} = T_{n}$. Hence there exists an element
$g_{0} \in \prod_{n} SL(d_{n}+1,q_{n})$ such that
$\prod_{n} T_{n} \leqslant G_{0} = \langle \prod_{n} S_{n} , g_{0} \rangle$.
Let $U_{n}$ be the subgroup of strictly upper triangular matrices
in $SL(d_{n}+1,q_{n})$, and let $H_{n}$ be the subgroup of diagonal
matrices. Then $B_{n} = U_{n} \rtimes H_{n}$ is a Borel subgroup of
$SL(d_{n}+1,q_{n})$. We shall show that
$\prod_{n} B_{n} \leqslant G_{0}$.

First we shall show that $\prod_{n} H_{n} \leqslant G_{0}$. Fix some
$n \in \mathbb{N}$. Let $D = \diag ( \lambda_{1}, \dots , 
\lambda_{d_{n}+1} ) \in H_{n}$. Then
$D_{1} = \diag ( \lambda_{1}, \lambda_{1}^{-1}, 1 , \dots , 1 )
\in H_{n} \cap S_{n}$,
$D_{2} = \diag ( 1, \lambda_{1} \lambda_{2}, \lambda_{3}, \dots ,
\lambda_{d_{n}+1} ) \in H_{n} \cap T_{n}$ and 
$D = D_{1} D_{2}$. The result follows.

Next we shall show that $\prod_{n} U_{n} \leqslant G_{0}$. Fix some
$n \in \mathbb{N}$. Note that
\[
\begin{pmatrix}
1 & \boldsymbol{0} & 0 \\
\boldsymbol{0} & A & \boldsymbol{b} \\
0 & \boldsymbol{0} & 1
\end{pmatrix}
\begin{pmatrix}
1 & \boldsymbol{c} & 0 \\
\boldsymbol{0} & \boldsymbol{1} & \boldsymbol{0} \\
0 & \boldsymbol{0} & 1
\end{pmatrix}
=
\begin{pmatrix}
1 & \boldsymbol{c} & 0 \\
\boldsymbol{0} & A & \boldsymbol{b} \\
0 & \boldsymbol{0} & 1
\end{pmatrix}
\]
for each $(d_{n}-1) \times (d_{n}-1)$-matrix $A$. Also note that if 
$Z \in U_{n}$ has the form
\[
\begin{pmatrix}
1 & \boldsymbol{0} & d \\
\boldsymbol{0} & \boldsymbol{1} & \boldsymbol{0} \\
0 & \boldsymbol{0} & 1
\end{pmatrix}
\]
then there exist $X \in U_{n} \cap S_{n}$ and $Y \in U_{n} \cap T_{n}$
such that $[X,Y] = Z$. (For example, this follows from Chevalley's
commutator formula \cite[5.2.2]{ca1}.) Hence if $\phi \in U_{n}$ is
arbitrary, then there exist $\theta$, $\tau \in U_{n} \cap S_{n}$ and
$\psi$, $\sigma \in U_{n} \cap T_{n}$ such that
$\phi = \psi \theta [ \tau, \sigma ]$. The result follows.

Let $N_{n}$ be the subgroup of $SL(d_{n}+1,q_{n})$ consisting of
the elements which stabilise
the frame $\{ \langle v_{i} \rangle \mid 1 \leq i \leq d_{n}+1 \}$.
Then the Weyl group of $SL(d_{n}+1,q_{n})$ is 
$W_{n} = N_{n}/B_{n} \cap N_{n}$; and $W_{n}$ is isomorphic to
$Sym(d_{n}+1)$ acting on the set
$\{ v_{i} \mid 1 \leq i \leq d_{n}+1 \}$. Note that $N_{n} \cap S_{n}$
corresponds to the subgroup $Sym(d_{n})$ of $W_{n}$. Let
$\theta = (\, d_{n} \, d_{n}+1 \,)$. Arguing as in the proof of Lemma
\ref{L:uni1}, we see that for every $\phi \in Sym(d_{n}+1)$, there exist
$\psi_{1}$, $\psi_{2}$, $\psi_{3} \in Sym(d_{n})$ such that
$\phi = \psi_{1} \theta \psi_{2} \theta \psi_{3}$. Hence
there exists $g_{1} \in \prod_{n} SL(d_{n}+1,q_{n})$ such that
$\prod_{n} N_{n} \leqslant G_{1} = \langle G_{0}, g_{1} \rangle$.
It follows that $G_{1} = \prod_{n} SL(d_{n}+1,q_{n})$. 
\end{proof}

Finally we shall prove the analogue of Lemma \ref{L:uni2}. Consider a
product of the form $\prod_{n} SL(8d_{n},q_{n})$. Fix some $n \in \mathbb{N}$.
Let $SL(8d_{n},q_{n})$ act on the vector space $V_{n} = V(8d_{n},q_{n})$
in the natural manner, and let $\mathcal{B}_{n} = \{ v_{i} \mid 1 \leq i \leq 8d_{n} \}$
be a basis of $V_{n}$. Let 
$E_{0} = \langle v_{i} \mid 1 \leq i \leq 4d_{n} \rangle$ and
$E_{1} = \langle v_{i} \mid 4d_{n}+1 \leq i \leq 8d_{n} \rangle$. We regard
$SL(E_{0})$ as the subgroup of $SL(8d_{n},q_{n})$ consisting of the elements
$\pi$ such that $\pi [ E_{0} ] = E_{0}$ and such that 
$\pi(v_{i}) = v_{i}$ for all $4d_{n}+1 \leq i \leq 8d_{n}$. We also
regard $SL(E_{1})$ as a subgroup of $SL(8d_{n},q_{n})$ in the obvious
fashion. Let $\Gamma_{n} = SL(E_{0}) \times SL(E_{1}) \leqslant
SL(8d_{n},q_{n})$.

\begin{Lem} \label{L:double}
$\prod_{n} SL(8d_{n},q_{n})$ is finitely generated over
$\prod_{n} \Gamma_{n}$.
\end{Lem}

\begin{proof}
Once again, we shall make use of the Bruhat decomposition of the special
linear group. Fix some $n \in \mathbb{N}$. We shall regard $SL(8d_{n},q_{n})$
as a group of matrices with respect to the ordered basis
$(v_{1}, \dots ,v_{8d_{n}})$ of $V_{n}$. Let $B_{n} = U_{n} \rtimes H_{n}$
be the Borel subgroup consisting of the upper triangular matrices of
$SL(8d_{n},q_{n})$ First we shall show that there exists a subgroup
$G_{0}$ of $\prod_{n} SL(8d_{n},q_{n})$ such that
\begin{enumerate}
\item $G_{0}$ is finitely generated over $\prod_{n} \Gamma_{n}$, and
\item $\prod_{n} U_{n} \leqslant G_{0}$.
\end{enumerate}
Fix some $n \in \mathbb{N}$. Let $M_{n}$ be the ring of all
$4d_{n} \times 4d_{n}$-matrices over $GF(q_{n})$, and let
\[
T_{n} = \left \{
\begin{pmatrix}
\boldsymbol{1} & S \\
\boldsymbol{0} & \boldsymbol{1}
\end{pmatrix}
\mid S \in M_{n} \right \}.
\]
Then it is enough to find $G_{0}$ such that $\prod_{n} T_{n} \leqslant G_{0}$.
Note that for each $A \in SL(4d_{n},q_{n})$, we have that
\[
\begin{pmatrix}
A & \boldsymbol{0} \\
\boldsymbol{0} & \boldsymbol{1}
\end{pmatrix}
\begin{pmatrix}
\boldsymbol{1} & S \\
\boldsymbol{0} & \boldsymbol{1}
\end{pmatrix}
\begin{pmatrix}
A^{-1} & \boldsymbol{0} \\
\boldsymbol{0} & \boldsymbol{1}
\end{pmatrix}
=
\begin{pmatrix}
\boldsymbol{1} & AS \\
\boldsymbol{0} & \boldsymbol{1}
\end{pmatrix}
.
\]
Regard $M_{n}$ as a $SL(4d_{n},q_{n})$-module with the natural
action, $S \overset{A}{\longmapsto} AS$. Then the existence of a
suitable subgroup $G_{0}$ is an immediate consequence of the
following claim.

\begin{Claim} \label{C:double}
$\prod_{n} M_{n}$ is finitely generated as a 
$\prod_{n} SL(4d_{n},q_{n})$-module.
\end{Claim}

\begin{proof}[Proof of Claim]
We shall prove that $M_{n}$ is ``uniformly generated'' as a
$SL(4d_{n},q_{n})$-module. The result will then follow. Fix
some $n \in \mathbb{N}$. Throughout this proof, each of the matrices
will be expressed in terms of $2d_{n} \times 2d_{n}$-blocks. Let
\[
J_{1} =
\begin{pmatrix}
\boldsymbol{1} & \boldsymbol{0} \\
\boldsymbol{0} & \boldsymbol{0}
\end{pmatrix}
\quad
J_{2} = 
\begin{pmatrix}
\boldsymbol{0} & \boldsymbol{0} \\
\boldsymbol{0} & \boldsymbol{1}
\end{pmatrix}
\quad
J_{3} = 
\begin{pmatrix}
\boldsymbol{0} & \boldsymbol{1} \\
\boldsymbol{0} & \boldsymbol{0}
\end{pmatrix}
\quad
J_{4} =
\begin{pmatrix}
\boldsymbol{0} & \boldsymbol{0} \\
\boldsymbol{1} & \boldsymbol{0}
\end{pmatrix}
.
\]
If $B \in GL(2d_{n},q_{n})$, then $
\begin{pmatrix}
B & \boldsymbol{0} \\
\boldsymbol{0} & B^{-1}
\end{pmatrix}
$, $
\begin{pmatrix}
B^{-1} & \boldsymbol{0} \\
\boldsymbol{0} & B
\end{pmatrix}
\in SL(4d_{n},q_{n})$; and
\[
\begin{pmatrix}
B & \boldsymbol{0} \\
\boldsymbol{0} & B^{-1}
\end{pmatrix}
\begin{pmatrix}
\boldsymbol{1} & \boldsymbol{0} \\
\boldsymbol{0} & \boldsymbol{0}
\end{pmatrix}
=
\begin{pmatrix}
B & \boldsymbol{0} \\
\boldsymbol{0} & \boldsymbol{0}
\end{pmatrix}
\text{ and }
\begin{pmatrix}
B^{-1} & \boldsymbol{0} \\
\boldsymbol{0} & B
\end{pmatrix}
\begin{pmatrix}
\boldsymbol{0} & \boldsymbol{0} \\
\boldsymbol{1} & \boldsymbol{0}
\end{pmatrix}
=
\begin{pmatrix}
\boldsymbol{0} & \boldsymbol{0} \\
B & \boldsymbol{0}
\end{pmatrix}
\text{ etc.}
\]
Thus if $B_{1}, \dots , B_{4} \in GL(2d_{n},q_{n})$, then there exist
$C_{1}, \dots , C_{4} \in SL(4d_{n},q_{n})$ such that
\[
\begin{pmatrix}
B_{1} & B_{2} \\
B_{3} & B_{4}
\end{pmatrix}
= \sum_{i=1}^{4} C_{i} J_{i} .
\]
Now suppose that $
\begin{pmatrix}
S_{1} & S_{2} \\
S_{3} & S_{4}
\end{pmatrix}
\in M_{n}$ is arbitrary. By \cite{ze}, each of the matrices $S_{i}$
is the sum of two non-singular ones. Hence there exist
$C_{1}, \dots , C_{4}, D_{1}, \dots ,D_{4} \in SL(4d_{n},q_{n})$
such that
\[
\begin{pmatrix}
S_{1} & S_{2} \\
S_{3} & S_{4}
\end{pmatrix}
= \sum_{i=1}^{4} C_{i} J_{i} + \sum_{i=1}^{4} D_{i} J_{i} .
\]
Thus each $M_{n}$ is ``uniformly generated'' from the generators
$J_{1}, \dots , J_{4}$.
\end{proof}

Next we shall show that there exists an element 
$g_{0} \in \prod_{n} SL(8d_{n},q_{n})$ such that
$\prod_{n} H_{n} \leqslant G_{1} = \langle G_{0},g_{0} \rangle$; and
hence $\prod_{n} B_{n} \leqslant G_{1}$. For each $\lambda \in
GF(q_{n})^{*}$, let $D_{\lambda} =
\diag(\lambda, 1, \dots, 1) \in GL(4d_{n},q_{n})$. Define
\[
F_{n} = \left \{
\begin{pmatrix}
D_{\lambda} & \boldsymbol{0} \\
\boldsymbol{0} & D_{\lambda}^{-1}
\end{pmatrix}
\mid \lambda \in GF(q_{n})^{*} \right \}.
\]
Since $\prod_{n} \Gamma_{n} \leqslant G_{0}$, it is enough to find
an element $g_{0}$ such that $\prod_{n} F_{n} \leqslant
\langle G_{0},g \rangle$. For each $\lambda \in GF(q_{n})^{*}$, let
$E_{\lambda} = \diag( \lambda, \lambda^{-1}, 1, \dots, 1) \in
GL(4d_{n},q_{n})$. Define
\[
K_{n} = \left \{
\begin{pmatrix}
E_{\lambda} & \boldsymbol{0} \\
\boldsymbol{0} & \boldsymbol{1}
\end{pmatrix}
\mid \lambda \in GF(q_{n})^{*} \right \}.
\]
Then $\prod_{n} K_{n} \leqslant \prod_{n} \Gamma_{n} \leqslant G_{0}$.
Also there exists an element $\pi \in N_{n}$ such that
$\pi K_{n} \pi^{-1} = F_{n}$. The existence of a suitable element
$g_{0}$ follows easily.

Finally we shall show that there exists an element 
$g_{1} \in \prod_{n} SL(8d_{n},q_{n})$ such that
$\prod_{n} N_{n} \leqslant G_{2} = \langle G_{1}, g_{1} \rangle$.
Fix some $n \in \mathbb{N}$. Let
$\mathcal{E}^{n}_{0} = \{ v_{i} \mid 1 \leq i \leq 4d_{n} \}$ and
$\mathcal{E}^{n}_{1} = \{ v_{i} \mid 4d_{n}+1 \leq i \leq 8d_{n} \}$;
so that $\mathcal{B}_{n} = \mathcal{E}^{n}_{0} \cup \mathcal{E}^{n}_{1}$.
Then the groups of permutation matrices $Alt( \mathcal{E}^{n}_{0})$,
$Alt(\mathcal{E}^{n}_{1})$ are subgroups of $\Gamma_{n}$. Lemma
\ref{L:uni2} implies that there exists an element
$g_{1} \in \prod_{n} Alt(\mathcal{B}_{n})$ such that
$\prod_{n} Alt(\mathcal{B}_{n}) = \langle 
\prod \left( Alt(\mathcal{E}^{n}_{0}) \times Alt(\mathcal{E}^{n}_{1}) \right),
g_{1} \rangle$. Now let $\pi = \langle \pi(n) \rangle_{n} \in
\prod_{n} N_{n}$ be an arbitrary element. Let $X$ be the subset of
$\mathbb{N}$ consisting of those $n$ such that $\pi(n)$ corresponds to
an odd permutation of $\mathcal{B}_{n}$.
Then there
exists an element $\psi_{X} = \langle \psi_{X}(n) \rangle_{n} \in
\prod_{n} \left( N_{n} \cap \Gamma_{n} \right)$ such that
\begin{enumerate}
\item $\psi_{X}(n)$ corresponds to the odd permutation
$(\, v_{1} \, v_{2} \,)$ if $n \in X$, and
\item $\psi_{X}(n) = 1$ if $n \notin X$.
\end{enumerate}
Since $\prod_{n}H_{n}$, $\prod_{n}Alt(\mathcal{B}_{n}) \leqslant
\langle G_{1}, g_{1} \rangle$, it follows that
$\pi \psi_{X} \in \langle G_{1}, g_{1} \rangle$; and hence
$\pi \in \langle G_{1}, g_{1} \rangle$.
Thus $\prod_{n} N_{n} \leqslant
\langle G_{1}, g_{1} \rangle$.
\end{proof}

Now we are ready to begin the proof of Theorem \ref{T:lin}. Since the
proof is very similar to that of Theorem \ref{T:alt}, we shall just
sketch the main points. Suppose that
$\langle SL(d_{n}, q_{n}) \mid n \in \mathbb{N} \rangle$ is a sequence
of finite special linear groups which satisfies the hypotheses of
Theorem \ref{T:lin}. By Proposition \ref{P:finite}, we can suppose
that $\{ d_{n} \mid n \in \mathbb{N} \}$ is an infinite subset of
$\mathbb{N}$. Arguing as in the proof of Lemma \ref{L:red}, we can
reduce to the case when each $d_{n}$ has the form $2^{m_{n}}$ for
some $m_{n} \geq 2$. (Since the sequence satisfies condition
\ref{T:lin}(2), there exists a finite set $\mathcal{F}$ of groups
such that if $d_{n} \leq 3$, then $SL(d_{n},q_{n}) \in \mathcal{F}$.
By Propositions \ref{P:gen} and \ref{P:finite}, we can safely
ignore these factors.) Let $G = 
\prod_{n} SL(2^{m_{n}},q_{n})$ and suppose that $c(G) = \omega$.
Express $G = \underset{t < \omega}{\bigcup} H_{t}$ as the union
of a chain of $\omega$ proper subgroups. Now suppose that
$t \geq 0$ and that we have defined
\begin{enumerate}
\item a sequence of elements
$f_{j} = \langle f_{j}(n) \rangle_{n} \in G$,
\item a strictly increasing sequence of integers $i_{j}$ such that
$f_{j} \in H_{i_{j}}$,
\item a sequence of elements
$g_{j} = \langle g_{j}(n) \rangle_{n} \in G \smallsetminus H_{i_{j}}$, and
\item a sequence of elements
$h_{j} = \langle h_{j}(n) \rangle_{n} \in G$
\end{enumerate}
for $0 \leq j \leq t$ such that the following conditions hold.
\begin{enumerate}
\item[(a)] $f_{0}(n), \dots , f_{t}(n)$ is a generic sequence in
$SL(2^{m_{n}},q_{n})$ for each $n \in \mathbb{N}$.
\item[(b)] If $0 \leq j \leq k \leq t$, then
$h_{k} f_{j} h_{k}^{-1} = g_{j}$.
\item[(c)] If $m_{n}-1 \leq j \leq t$, then
$h_{j}(n) = h_{m_{n}-1}(n)$.
\end{enumerate}
We must show that it is possible to continue the induction. There is
no difficulty in defining $f_{t+1}$ and $i_{t+1}$. The problem is to show
that there exist suitable elements $h_{t+1}$ and
$g_{t+1} = h_{t+1} f_{t+1} h_{t+1}^{-1}$ such that 
$g_{t+1} \notin H_{i_{t+1}}$. As in the proof of Theorem \ref{T:alt},
we shall show that if no such elements exist, then $G$ is finitely
generated over $H_{i_{t+1}}$; which is a contradiction. So suppose
that no such elements exist. Let
$P^{*} = \prod \{ SL(2^{m_{n}},q_{n}) \mid m_{n} \leq t+3 \}$ and
$P^{\prime} = \prod \{ SL(2^{m_{n}},q_{n}) \mid m_{n} \geq t+4 \}$;
so that $G = P^{*} \times P^{\prime}$.
Since $\langle SL(2^{m_{n}},q_{n} \mid n \in \mathbb{N} \rangle$ satisfies
condition \ref{T:lin}(2), either $c(P^{*}) > \omega$ or $P^{*}$ is
finite. Thus we can suppose that $i_{t+1}$ was chosen so that
$P^{*} \leqslant H_{i_{t+1}}$. Fix some $n \in \mathbb{N}$ such that
$m_{n} \geq t+4$. Let 
$\mathcal{B}_{n}$ be a basis
of $V(2^{m_{n}},q_{n})$ chosen so that the group $E(n) =
\langle g_{0}(n), \dots , g_{t}(n) \rangle$ is an elementary abelian
subgroup of $Alt(\mathcal{B}_{n})$ of order $2^{t+1}$, 
which acts semiregularly on $\mathcal{B}_{n}$. Let
$\{ \Phi^{n}_{i} \mid 1 \leq i \leq 2^{m_{n}-t-1} \}$ be the set of
orbits of $E(n)$ on $\mathcal{B}_{n}$. For each
$1 \leq i \leq 2^{m_{n}-t-1}$, choose $v^{1}_{i} \in \Phi^{n}_{i}$.
Let $E(n) = \{ \psi_{k} \mid 1 \leq k \leq 2^{t+1} \}$. For each
$1 \leq k \leq 2^{t+1}$ and $1 \leq i \leq 2^{m_{n}-t-1}$, define
$v^{k}_{i} = \psi_{k}(v^{1}_{i})$. For each $1 \leq k \leq 2^{t+1}$,
let $V^{n}_{k} = \langle v^{k}_{i} \mid 1 \leq i \leq
2^{m_{n}-t-1} \rangle$. Then $V(2^{m_{n}},q_{n}) = V^{n}_{1} \oplus \dots
\oplus V^{n}_{2^{t+1}}$, and the diagonal subgroup
\[
D_{n} = Diag \left( SL(V^{n}_{1}) \times \dots \times SL(V^{n}_{2^{t+1}})
\right)
\]
is contained in the centraliser of
$\langle g_{0}(n), \dots , g_{t}(n) \rangle$ in $SL(2^{m_{n}},q_{n})$.
Since each candidate $\pi$ for $g_{t+1}$ satisfies $\pi \in H_{i_{t+1}}$,
we find that the following statement holds.
\begin{enumerate}
\item[(\dag)] Suppose that $\pi \in P^{\prime}$. If
$\pi(n) \in D_{n} \cap C(2^{m_{n}},q_{n})$ for all $n$ such that
$m_{n} \geq t+4$, then $\pi \in H_{i_{t+1}}$.
\end{enumerate}
Using Theorem \ref{T:saxl}, this implies that
$\prod \{ D_{n} \mid m_{n} \geq t+4 \} \leqslant H_{i_{t+1}}$.
Arguing as in the proof of Theorem \ref{T:alt}, we see that there
exists a subgroup $\Gamma_{0}$ of $G$ such that
\begin{enumerate}
\item $\Gamma_{0}$ is finitely generated over $H_{i_{t+1}}$, and
\item $\prod \{ SL(V^{n}_{1}) \times \dots \times SL(V^{n}_{2^{t+1}})
\mid m-{n} \geq t+4 \} \leqslant \Gamma_{0}$.
\end{enumerate}
By repeatedly applying Lemma \ref{L:double}, we next see that there
exists a subgroup $\Gamma_{1}$ of $G$ such that
\begin{enumerate}
\item $\Gamma_{1}$ is finitely generated over $\Gamma_{0}$, and
\item $P^{\prime} = \prod \{ SL(2^{m_{n}},q_{n}) \mid m_{n} \geq t+4 \}
\leqslant \Gamma_{1}$.
\end{enumerate}
But this means that $\Gamma_{1} = G$; and so $G$ is finitely generated
over $H_{i_{t+1}}$. This contradiction shows that the induction can be 
carried out for all $t < \omega$. But this yields an element $h \in G$
such that $h f_{t} h^{-1} = g_{t}$ for all $t < \omega$, which is
impossible. Thus $c(G) > \omega$. This completes the proof of
Theorem \ref{T:lin}.

\section{The proof of Theorem \ref{T:class}} \label{S:class}
In this section, we shall complete the proof of Theorem \ref{T:class}.
Most of our work will go into proving the special cases of
Theorem \ref{T:class} in which each $S_{n}$ is a classical group of a
fixed kind. We shall deal successively with the symplectic groups,
the unitary groups and the orthogonal groups over finite fields.
The general result will then follow easily. (Clear accounts of the
classical groups can be found in \cite{ca1} and \cite{ta}.)

\subsection{Symplectic groups} \label{SS:symp}
Suppose that $\langle Sp(2d_{n},q_{n}) \mid n \in \mathbb{N} \rangle$
is a sequence of finite symplectic groups such that $d_{n} \geq 2$
for each $n \in \mathbb{N}$. 
Fix some $n \in \mathbb{N}$.
Then there exists a basis
$\boldsymbol{e}\sphat \, \boldsymbol{f}
= (e_{i} \mid 1 \leq i \leq d_{n})\sphat \, (f_{i} \mid 1 \leq i \leq d_{n})$
of the corresponding symplectic space such that
$(e_{i},f_{j}) = \delta_{ij}$ and
$(e_{i},e_{j}) = (f_{i},f_{j}) =0$ for all
$1 \leq i,j \leq d_{n}$. (Such a basis is called a {\em normal basis\/}.)
We shall consider $Sp(2d_{n},q_{n})$ as a group of matrices with respect
to the ordered basis $\boldsymbol{e} \sphat \, \boldsymbol{f}$. Let
$E_{d_{n}} = \langle e_{1}, \dots , e_{d_{n}} \rangle$ and
$F_{d_{n}} = \langle f_{1}, \dots , f_{d_{n}} \rangle$. Then the setwise
stabiliser of the subspaces $E_{d_{n}}$ and $F_{d_{n}}$ in $Sp(2d_{n},q_{n})$
contains the subgroup
\[
G_{n} = \left\{
\begin{pmatrix}
A & \boldsymbol{0} \\
\boldsymbol{0} & \left( A^{-1} \right)^{T}
\end{pmatrix}
\mid A \in SL(d_{n},q_{n}) \right\}.
\]

\begin{Thm} \label{T:symp}
Suppose that $d_{n} \geq 3$ for all $n \in \mathbb{N}$. Then
$\prod_{n} Sp(2d_{n},q_{n})$ is finitely generated over the
subgroup $\prod_{n} G_{n}$.
\end{Thm}

\begin{Cor} \label{C:symp}
Suppose that $\langle S_{n} \mid n \in \mathbb{N} \rangle$ is a sequence
of finite simple symplectic groups such that there does {\em not\/}
exist an infinite subset $I$ of $\mathbb{N}$ for which conditions
\ref{T:count}(1) and \ref{T:count}(2) are satisfied. Then
$c(\prod_{n} S_{n}) > \omega$.
\end{Cor}

\begin{proof}[Proof of Corollary \ref{C:symp}]
For each $n \in \mathbb{N}$, let $S_{n} = PSp(2d_{n},q_{n})$. Put
$J = \{n \in \mathbb{N} \mid d_{n} < 3 \}$, so that
$\prod_{n \in \mathbb{N}} S_{n} = 
\left( \prod_{n \in J} S_{n} \right) \times
\left( \prod_{n \notin J} S_{n} \right)$.
By assumption, there exists a finite set $\mathcal{F}$ of groups such that
$S_{n} \in \mathcal{F}$ for all $n \in J$.
By Proposition \ref{P:finite},
either $c(\prod_{n \in J} S_{n}) > \omega$ or $\prod_{n \in J} S_{n}$
is finite. Hence if $\prod_{n \notin J} S_{n}$ is finite, then the
result follows from Proposition \ref{P:gen}. So we can suppose that
$\prod_{n \notin J} S_{n}$ is infinite; and 
it is enough to prove
that $c(\prod_{n \notin J} S_{n}) > \omega$. To simplify
notation, we shall suppose that $J = \emptyset$. Let 
$\prod_{n} G_{n}$ be the subgroup of $\prod_{n} Sp(2d_{n},q_{n})$
defined above. By Theorem \ref{T:lin}, $c(\prod_{n} G_{n})
> \omega$. So using Theorem \ref{T:symp} and Proposition \ref{P:gen},
we see that $c(\prod_{n} Sp(2d_{n},q_{n})) > \omega$. Hence
$c(\prod_{n} PSp(2d_{n},q_{n}) ) > \omega$.
\end{proof}

We shall approach Theorem \ref{T:symp} via the Bruhat decomposition
\[
Sp(2d,q) = \underset{w \in W}{\bigcup} BwB
\]
of the symplectic group, where $B$ is a Borel subgroup and $W$ is the
Weyl group. Fix some $n \in \mathbb{N}$. Let
$\boldsymbol{e} \sphat \, \boldsymbol{f} = 
(e_{i} \mid 1 \leq i \leq d_{n}) \sphat \,
(f_{i} \mid 1 \leq i \leq d_{n})$ be our distinguished normal basis.
For each $1 \leq i \leq d_{n}$, let
$E_{i} = \langle e_{1}, \dots e_{i} \rangle$. Then the stabiliser
$B_{n}$ of the flag of totally isotropic subspaces
\[
E_{1} \leqslant E_{2} \leqslant \dots \leqslant E_{d_{n}}
\]
is a Borel subgroup of $Sp(2d_{n},q_{n})$. Let $N_{n}$ be the 
subgroup of $Sp(2d_{n},q_{n})$ which stabilises the symplectic
frame $\{ \langle e_{i} \rangle, \langle f_{i} \rangle \mid
1 \leq i \leq d_{n} \}$. Then the Weyl group of $Sp(2d_{n},q_{n})$
is $N_{n}/B_{n} \cap N_{n}$. Let $H_{n} = B_{n} \cap N_{n}$.
Then $H_{n}$ consists of the matrices of the form
\[
\begin{pmatrix}
D & \boldsymbol{0} \\
\boldsymbol{0} & D^{-1}
\end{pmatrix}
\]
where $D \in GL(d_{n},q_{n})$ is a diagonal matrix. 
Let $UT_{n}$ be the subgroup of strictly
upper triangular matrices in $SL(d_{n},q_{n})$, and define
\[
U_{n} = \left\{
\begin{pmatrix}
P & PS \\
\boldsymbol{0} & \left( P^{-1} \right)^{T}
\end{pmatrix}
\mid P \in UT_{n}, S^{T} = S \right\}.
\]
Then $B_{n} = U_{n} \rtimes H_{n}$.
 
First we shall show that there exists a subgroup $\Gamma_{0}$ of
$\prod_{n} Sp(2d_{n},q_{n})$ such that
\begin{enumerate}
\item $\Gamma_{0}$ is finitely generated over $\prod_{n} G_{n}$, and
\item $\prod_{n} U_{n} \leqslant \Gamma_{0}$.
\end{enumerate}
Note that for each
$A \in SL(d_{n},q_{n})$, we have that
\[
\begin{pmatrix}
A & \boldsymbol{0} \\
\boldsymbol{0} & \left( A^{-1} \right)^{T}
\end{pmatrix}
\begin{pmatrix}
\boldsymbol{1} & S \\
\boldsymbol{0} & \boldsymbol{1}
\end{pmatrix}
\begin{pmatrix}
A^{-1} & \boldsymbol{0} \\
\boldsymbol{0} & A^{T}
\end{pmatrix}
=
\begin{pmatrix}
\boldsymbol{1} & ASA^{T} \\
\boldsymbol{0} & \boldsymbol{1}
\end{pmatrix}
.
\]
Let $M_{n}$ be the left $SL(d_{n},q_{n})$-module of symmetric
$d_{n} \times d_{n}$-matrices, with the action
\[
S \overset{A}{\longmapsto} ASA^{T} .
\]
Then it is enough to prove that $\prod_{n} M_{n}$ is finitely
generated as a $\prod_{n} SL(d_{n},q_{n})$-module. We shall consider
$M_{n}$ in three different cases, and show that in each case $M_{n}$
is ``uniformly generated'' as a $SL(d_{n},q_{n})$-module. The result
will then follow. Let $p = \operatorname{char} (GF(q_{n}))$.

\noindent
{\em Case 1.\/} Suppose that $p > 3$. Since $p$ is odd, every $S \in M_{n}$
is congruent to a diagonal matrix. This easily implies that there exists
$A \in SL(d_{n},q_{n})$ such that $ASA^{T}$ is a diagonal matrix. Thus
we need only consider diagonal matrices
$D = \diag (\lambda_{1}, \dots , \lambda_{d_{n}} ) \in M_{n}$. Let
$D_{1} = \diag (1, 0, \dots , 0 )$ and
$D_{2} = \diag (0, 1, \dots , 1 )$, so that
$\boldsymbol{1} = D_{1} + D_{2}$. If
$\boldsymbol{\alpha} = ( \alpha_{1}, \dots , \alpha_{d_{n}})
\in \left( GF(q_{n})^{*} \right)^{d_{n}}$, let
$R_{\boldsymbol{\alpha}} = \diag (\alpha_{1}, 1, \dots , 1, \alpha^{-1}_{1})$
and $S_{\boldsymbol{\alpha}} =
\diag ( (\alpha_{2} \dots \alpha_{d_{n}})^{-1}, \alpha_{2}, \dots ,
\alpha_{d_{n}})$. Then
$R_{\boldsymbol{\alpha}} D_{1} R^{T}_{\boldsymbol{\alpha}} =
\diag ( \alpha^{2}_{1}, 0, \dots , 0)$ and
$S_{\boldsymbol{\alpha}} D_{2} S_{\boldsymbol{\alpha}}^{T} =
\diag (0, \alpha_{2}^{2}, \dots , \alpha_{d_{n}}^{2} )$. Since
$p > 3$, for each $\lambda \in GF(q_{n})$, there exist
$\beta_{1}, \dots , \beta_{4} \in GF(q_{n})^{*}$ such that
$\lambda = \sum_{i = 1}^{4} \beta_{i}^{2}$. (For example, see
Chapter 4 \cite{sm}.) Thus we can ``uniformly generate'' each
diagonal matrix $D \in M_{n}$ from the generators $D_{1}$ and
$D_{2}$.

\noindent
{\em Case 2.\/} Suppose that $p = 3$. Once again, we need only consider
diagonal matrices $D \in M_{n}$.
Let $D_{1} = \diag(1, 1, 0, \dots , 0)$ and 
$D_{2} = \diag(\delta, 1, 1, \dots , 1)$, where $\delta = 1$ if
$d_{n}$ is even and $\delta = 0$ if $d_{n}$ is odd. For each subset
$X$ of $\{1, \dots , d_{n} \}$, let
$D_{1}^{X} = \diag( \chi_{X}(1), 0, \dots , 0)$ and
$D_{2}^{X} = \diag(0, \chi_{X}(2), \dots ,\chi_{X}(d_{n}) )$,
where $\chi_{X}$ is the characteristic function of $X$.
It is easy to check that if
\[
S \in \left\{
\begin{pmatrix}
0 & 0 \\
0 & 0
\end{pmatrix}
,
\begin{pmatrix}
1 & 0 \\
0 & 0 
\end{pmatrix}
,
\begin{pmatrix}
0 & 0 \\
0 & 1
\end{pmatrix}
,
\begin{pmatrix}
1 & 0 \\
0 & 1
\end{pmatrix}
\right\}
\]
then there exist $A_{1}, \dots , A_{6} \in SL(2,3)$ such that
$S = \sum_{i = 1}^{6} A_{i} A_{i}^{T}$. Hence there exist
$B_{i}$, $C_{i} \in SL(d_{n},3)$ for $1 \leq i \leq 6$ such that
$\sum_{i=1}^{6} B_{i} D_{1} B_{i}^{T} = D_{1}^{X}$ and
$\sum_{i=1}^{6} C_{i} D_{2} C_{i}^{T} = D_{2}^{X}$.  Now it is easy
to complete the proof of this case. Let 
$D = \diag( \lambda_{1}, \dots , \lambda_{d_{n}}) \in M_{n}$.
Then for each $1 \leq i \leq d_{n}$, there exist 
$\alpha_{i}$, $\beta_{i} \in GF(q_{n})$ such that
$\lambda_{i} = \alpha_{i}^{2} + \beta_{i}^{2}$. (It is well-known
that if $\mathbb{F}$ is {\em any\/} finite field, then every element
of $\mathbb{F}$ is a sum of two squares. Unfortunately it is often
not possible to express an element as the sum of two nonzero
squares.) Let $X = \{ i \mid \alpha_{i} \ne 0 \}$ and
$Y = \{ i \mid \beta_{i} \ne 0 \}$. Then there exist diagonal
matrices $R_{j} \in SL(d_{n}, q_{n})$ for $1 \leq j \leq 4$
such that
$R_{1} D_{1}^{X} R_{1}^{T} = \diag(\alpha_{1}^{2}, 0, \dots , 0)$,
$R_{2} D_{2}^{X} R_{2}^{T} = \diag( 0, \alpha_{2}^{2}, \dots ,
\alpha_{d_{n}}^{2} )$, $R_{3} D_{1}^{Y} R_{3}^{T} =
\diag (\beta_{1}^{2}, 0, \dots , 0)$ and $R_{4} D_{2}^{Y} R_{4}^{T} =
\diag(0, \beta_{2}^{2}, \dots , \beta_{d_{n}}^{2})$. Thus we can
``uniformly generate'' each diagonal matrix from the generators
$D_{1}$ and $D_{2}$.

\noindent
{\em Case 3.\/} Suppose that $p = 2$. If $S = (s_{ij}) \in M_{n}$,
then $S$ is said to be {\em alternating\/} if $s_{ii} = 0$ for
all $1 \leq i \leq d_{n}$. If $S$ is {\em not\/} alternating, then
$S$ is congruent to a diagonal matrix. Clearly for each $S \in M_{n}$,
there exist $B_{k} \in M_{n}$ for $1 \leq k \leq 3$ such that
$S = \sum_{k=1}^{3} B_{k}$ and none of the $B_{k}$ are alternating.
Thus we need only consider diagonal matrices $D \in M_{n}$. It is
easily checked that if $C$ is a diagonal $3 \times 3$-matrix over
$GF(2)$, then there exist $A_{i} \in SL(3,2)$ for
$1 \leq i \leq 4$ such that
$C = \sum_{i=1}^{4} A_{i} A_{i}^{T}$. It is now easy to adapt the
argument of Case 2. (In fact, the argument is even simpler in this case,
as every element in $GF(q_{n})$ is a square.)
This completes the proof of the existence of the subgroup $\Gamma_{0}$ of
$\prod_{n} Sp(2d_{n},q_{n})$.

Next we shall show that there exists a subgroup $\Gamma_{1}$ of
$\prod_{n} Sp(2d_{n},q_{n})$ such that
\begin{enumerate}
\item $\Gamma_{1}$ is finitely generated over $\Gamma_{0}$, and
\item $\prod_{n} H_{n} \leqslant \Gamma_{1}$.
\end{enumerate}
For each $n \in \mathbb{N}$, let $LT_{n}$ be the subgroup of strictly
lower triangular matrices in $SL(d_{n},q_{n})$, and define
\[
V_{n} = \left \{
\begin{pmatrix}
Q & \boldsymbol{0} \\
SQ & \left( Q^{-1} \right)^{T}
\end{pmatrix}
\mid Q \in LT_{n}, S^{T} = S \right \}.
\]
Then there exists an element $\pi \in Sp(2d_{n},q_{n})$ such that
$\pi U_{n} \pi^{-1} = V_{n}$. Hence there exists
$g_{1} \in \prod_{n} Sp(2d_{n},q_{n})$ such that 
$\prod_{n} V_{n} \leqslant \Gamma_{1} =
\langle \Gamma_{0}, g_{1} \rangle$. We shall prove that
$\prod_{n} H_{n} \leqslant \Gamma_{1}$. For each
$\lambda \in GF(q_{n})^{*}$, let
$D_{\lambda} = \diag (\lambda, 1, \dots , 1) \in GL(d_{n},q_{n})$.
Define
\[
F_{n} = \left \{
\begin{pmatrix}
D_{\lambda} & \boldsymbol{0} \\
\boldsymbol{0} & D_{\lambda}^{-1}
\end{pmatrix}
\mid \lambda \in GF(q_{n})^{*} \right \}.
\]
Since $\prod_{n} G_{n} \leqslant \Gamma_{1}$, it suffices to show that
$\prod_{n} F_{n} \leqslant \Gamma_{1}$. For each $t \in GF(q_{n})$, let
$S(t)$ be the symmetric $d_{n} \times d_{n}$-matrix with $t$ in the
upper left position and 0 elsewhere. Define
\[
X(t) =
\begin{pmatrix}
\boldsymbol{1} & S(t) \\
\boldsymbol{0} & \boldsymbol{1}
\end{pmatrix}
\in U_{n} \text{ and } T(t) =
\begin{pmatrix}
\boldsymbol{1} & \boldsymbol{0} \\
S(t) & \boldsymbol{1}
\end{pmatrix}
\in V_{n}.
\]
Then it is easily checked that for each $\lambda \in GF(q_{n})^{*}$,
\[
\begin{pmatrix}
D_{\lambda} & \boldsymbol{0} \\
\boldsymbol{0} & D_{\lambda}^{-1}
\end{pmatrix}
= X(\lambda)Y(-\lambda^{-1})X(\lambda)X(-1)Y(1)X(-1).
\]
( This is essentially \cite[6.4.4]{ca1}.) Hence we have that
$\prod_{n} F_{n} \leqslant \Gamma_{1}$, as required.

Finally we shall show that there exists a
subgroup $\Gamma_{2}$ of $\prod_{n} Sp(2d_{n},q_{n})$ such that
\begin{enumerate}
\item $\Gamma_{2}$ is finitely generated over $\Gamma_{1}$; and
\item $\prod_{n} N_{n} \leqslant \Gamma_{2}$.
\end{enumerate}
This implies that $\Gamma_{2} = \prod_{n} Sp(2d_{n},q_{n})$, and hence
completes the proof of Theorem \ref{T:symp}. Again fix some
$n \in \mathbb{N}$. Then $W_{n}$ is generated by the images of the
elements $\{ w_{i} \mid 1 \leq i \leq d_{n} \}$ of $N_{n}$ defined
as follows.
\begin{enumerate}
\item[(a)] If $1 \leq i < d_{n}$, then $w_{i}$ is the permutation
matrix corresponding to the permutation
$(\, e_{i} \, e_{i+1} \,)(\, f_{i} \, f_{i+1} \,)$.
\item[(b)] $w_{d_{n}}(e_{d_{n}}) = - f_{d_{n}}$,
$w_{d_{n}}(f_{d_{n}}) = e_{d_{n}}$ and $w_{d_{n}}$ fixes the 
remaining elements of
$\boldsymbol{e}\sphat\, \boldsymbol{f}$. (Thus $w_{d_{n}}$ 
corresponds to the odd permutation $(\, e_{d_{n}} \, f_{d_{n}} \,)$.)
\end{enumerate}
It follows that $W_{n}$ is isomorphic to
$\mathbb{Z}^{d_{n}}_{2} \rtimes Sym(d_{n})$, where $\mathbb{Z}^{d_{n}}_{2}$
is the natural permutation module for $Sym(d_{n})$. Let
$t = \left \lfloor d_{n} /2 \right \rfloor$ and let
$v \in \mathbb{Z}^{d_{n}}_{2}$ be a vector of weight $t$.
Let $E_{n}$ be the submodule of $\mathbb{Z}^{d_{n}}_{2}$ consisting
of the vectors of even weight. Then for every 
$u \in E_{n}$, there exist $\pi$, $\phi \in Sym(d_{n})$
such that $u = \pi(v) + \phi(v)$.

Let $W^{+}_{n}$ be the subgroup of $W_{n}$ consisting of the even
permutations of the set
$\{ e_{i} , f_{i} \mid 1 \leq i \leq d_{n} \}$.
Then $W^{+}_{n}$ can be regarded as a subgroup of $Sp(2d_{n},q_{n})$.
Also notice that $W^{+}_{n}$ corresponds to the subgroup
$E_{n} \rtimes Sym(d_{n})$ of $\mathbb{Z}^{d_{n}}_{2}
\rtimes Sym(d_{n})$. So the 
argument of the previous paragraph shows that there exists
$g_{2} \in \prod_{n} Sp(2d_{n},q_{n})$ such that
$\prod_{n} W^{+}_{n} \leqslant \langle \Gamma_{1}, g_{2} \rangle$. We shall
show that $\prod_{n} N_{n} \leqslant \langle \Gamma_{1},g_{2} \rangle$.
Once again fix some $n \in \mathbb{N}$. Consider the element
$w = w_{d_{n}} w_{d_{n}-1} w_{d_{n}}$. Since $w$ corresponds to
an even permutation of $\{ e_{i},f_{i} \mid 1 \leq i \leq d_{n} \}$,
it follows that $w \in \langle W^{+}_{n}, B_{n} \rangle$. Using
the standard properties of groups with $BN$-pairs \cite[Section 8.2]{ca1},
we have that
\[
\left( B_{n} w_{d_{n}} B_{n} \right) \left( B_{n} w B_{n} \right)
= B_{n} w_{d_{n}} w B_{n} \cup B_{n} w B_{n}.
\]
In particular, $w \in \left(B_{n}w_{d_{n}}B_{n}\right)
\left(B_{n}wB_{n}\right)$.
Hence there exist elements $b_{1}$, $b_{2}$, $b_{3} \in B_{n}$ such
that $w = b_{1} w_{d_{n}} b_{2} w b_{3}$, and so
$w_{d_{n}} = b_{1}^{-1} w b_{3}^{-1} w^{-1} b_{2}^{-1}$. 
Obviously if we let
$1=c_{1}=c_{2}=c_{3} \in B_{n}$, then
$1 = c_{1}^{-1} w c_{3}^{-1} w^{-1} c_{2}^{-1}$. Hence for each subset $X$ of
$\mathbb{N}$, $\psi_{X} = \langle \psi_{X}(n) \rangle_{n} \in 
\langle \Gamma_{1}, g_{2} \rangle$, where $\psi_{X}(n) = w_{d_{n}}$ if
$n \in X$ and $\psi_{X}(n) = 1$ if $n \notin X$. It follows
easily that $\prod_{n} N_{n} \leqslant \langle \Gamma_{1}, g_{2} \rangle$.
( Cf. the final paragraph of the proof of Lemma \ref{L:double}.)

\subsection{Unitary groups} \label{SS:unit}
In this subsection, we shall consider products 
$\prod_{n} SU(d_{n},q_{n})$ of finite special unitary groups. First consider
the case when $d_{n}$ is even. Then the corresponding unitary space
has a normal basis $\boldsymbol{e} \sphat\, \boldsymbol{f}$.
Arguing as in Subsection \ref{SS:symp}, we obtain the following
result.

\begin{Thm} \label{T:unit1}
Suppose that $\langle SU(2d_{n},q_{n}) \mid n \in \mathbb{N} \rangle$
is a sequence of special unitary groups which satisfies the following
conditions.
\begin{enumerate}
\item[(1)] If $d_{n} = 1$, then $q_{n} > 3$.
\item[(2)] There does {\em not\/} exist an infinite subset $I$ of $\mathbb{N}$
and an integer $d$ such that
\begin{enumerate}
\item $d_{n} = d$ for all $n \in I$; and
\item if $n$, $m \in I$ and $n < m$, then $q_{n} < q_{m}$.
\end{enumerate}
\end{enumerate}
Then $c(\prod_{n} SU(2d_{n},q_{n})) > \omega$.
\end{Thm}
\begin{flushright}
$\square$
\end{flushright}

Next consider a product of the form $\prod_{n} SU(2d_{n}+1,q_{n})$,
where $d_{n} \geq 2$ for all $n \in \mathbb{N}$. Fix some $n \in \mathbb{N}$.
Then there exists a basis
\[
\boldsymbol{e} \sphat \, \boldsymbol{f} \sphat \, (w) =
(e_{i} \mid 1 \leq i \leq d_{n} ) \sphat \,
(f_{i} \mid 1 \leq i \leq d_{n} ) \sphat \, (w)
\]
of the corresponding unitary space such that
\[
(e_{i}, f_{j}) = \delta_{ij} \text{ and }
(e_{i}, e_{j}) = (f_{i},f_{j}) = 0
\]
for all $1 \leq i,j \leq d_{n}$; and
\[
(w,w) = 1 \text{ and } (w,e_{i}) = (w,f_{i}) = 0
\]
for all $1 \leq i \leq d_{n}$. Then we can regard $SU(2d_{n},q_{n})$
as the subgroup of \\ 
$SU(2d_{n}+1,q_{n})$ consisting of the elements
$\pi$ such that $\pi(w) = w$.

\begin{Thm} \label{T:unit2}
Suppose that $d_{n} \geq 2$ for all $n \in \mathbb{N}$. Then 
$\prod_{n} SU(2d_{n}+1,q_{n})$ is finitely generated over
$\prod_{n} SU(2d_{n},q_{n})$.
\end{Thm}

\begin{proof}
We shall make use of the Bruhat decomposition
\[
SU(2d+1,q) = \underset{ w \in W}{\bigcup} BwB
\]
of the special unitary group, where $B$ is a Borel subgroup and $W$ is
the Weyl group. Fix some $n \in \mathbb{N}$. For each $1 \leq i \leq d_{n}$,
let $E_{i} = \langle e_{1}, \dots , e_{i} \rangle$. Then the stabiliser
$B_{n}$ of the flag of totally isotropic subspaces
\[
E_{1} \leqslant E_{2} \leqslant \dots \leqslant E_{d_{n}}
\]
is a Borel subgroup of $SU(2d_{n}+1,q_{n})$. Let $N_{n}$ be the subgroup of \\
$SU(2d_{n}+1,q_{n})$ which stabilises the polar frame
$\{ \langle e_{i} \rangle , \langle f_{i} \rangle \mid
1 \leq i \leq d_{n} \}$. Then the Weyl group of $SU(2d_{n}+1,q_{n})$
is $W_{n} = N_{n}/ B_{n} \cap N_{n}$. Note that
$N_{n} \cap SU(2d_{n},q_{n})$ already contains representatives of
each element of $W_{n}$. Thus it suffices to prove that there exists
an element $g$ such that 
$\prod_{n} B_{n} \leqslant \langle \prod_{n} SU(2d_{n},q_{n}), g
\rangle$. We shall regard $SU(2d_{n}+1,q_{n})$ as a group of 
matrices with respect to the ordered basis
$(e_{1}, \dots , e_{d_{n}}, w, f_{d_{n}}, \dots , f_{1})$.
Thus $B_{n}$ is the subgroup of the upper triangular matrices
which are contained in $SU(2d_{n}+1,q_{n})$. Let $U_{n}$ be the 
subgroup of $B_{n}$ consisting of the strictly upper triangular
matrices. Let $H_{n} = B_{n} \cap N_{n}$. Then $H_{n}$ consists
of the diagonal matrices of the form
\[
\begin{pmatrix}
D & \boldsymbol{0} & \boldsymbol{0} \\
\boldsymbol{0} & \lambda^{-1} & \boldsymbol{0} \\
\boldsymbol{0} & \boldsymbol{0} & D^{*}
\end{pmatrix}
\]
where
\begin{enumerate}
\item $D = \diag (\lambda_{1}, \dots , \lambda_{d_{n}} ) \in
GL(d_{n}, q_{n}^{2})$;
\item $D^{*} = \diag ( \bar{\lambda}_{d_{n}}^{-1}, \dots ,
\bar{\lambda}_{1}^{-1} )$; and
\item $\lambda = \det (D D^{*})$.
\end{enumerate}
Here $\sigma \longmapsto \bar{\sigma}$ is the automorphism of
$GF(q_{n}^{2})$ of order 2. (Notice that for all diagonal matrices
$D \in GL(d_{n}, q_{n}^{2})$, we have that $\lambda \bar{\lambda} = 1$
and hence $( \lambda^{-1}w, \lambda^{-1}w) = (w,w)$.) We have that
$B_{n} = U_{n} \rtimes H_{n}$.

First we shall show that there exists an element $g \in \prod_{n}
SU(2d_{n}+1,q_{n})$ such that
$\prod_{n} U_{n} \leqslant \langle \prod_{n} SU(2d_{n},q_{n}), g
\rangle$. Later we shall see that we also have that
$\prod_{n} H_{n} \leqslant \langle \prod_{n} SU(2d_{n},q_{n}),g
\rangle$; and so $g$ satisfies our requirements. Once more, fix
some $n \in \mathbb{N}$. Let 
$E_{d_{n}}^{+} = \langle E_{d_{n}}, w \rangle$, and let
$\Gamma_{n}$ be the setwise stabiliser of $E_{d_{n}}^{+}$
in $SU(2d_{n}+1,q_{n})$. Let
$\rho : \Gamma_{n} \to GL( E_{d_{n}}^{+})$ be the restriction
map. We shall regard $GL( E_{d_{n}}^{+})$ as a group of matrices
with respect to the ordered basis
$( e_{1}, \dots ,e_{d_{n}},w)$. Note that for each
$A \in SL(d_{n}, q_{n}^{2})$, we have that 
\[
\begin{pmatrix}
A & \boldsymbol{0} \\
\boldsymbol{0} & 1
\end{pmatrix}
\in \rho \left[ SU(2d_{n},q_{n}) \cap \Gamma_{n} \right] ;
\]
and that for each $\boldsymbol{x} \in GF(q_{n}^{2})^{d_{n}}$,
we have that
\[
\begin{pmatrix}
A & \boldsymbol{0} \\
\boldsymbol{0} & 1
\end{pmatrix}
\begin{pmatrix}
\boldsymbol{1} & \boldsymbol{x} \\
\boldsymbol{0} & 1
\end{pmatrix}
\begin{pmatrix}
A^{-1} & \boldsymbol{0} \\
\boldsymbol{0} & 1
\end{pmatrix}
=
\begin{pmatrix}
\boldsymbol{1} & A \boldsymbol{x} \\
\boldsymbol{0} & 1
\end{pmatrix}
.
\]
Now choose $\theta \in U_{n}$ such that
$\rho ( \theta ) =
\begin{pmatrix}
\boldsymbol{1} & \boldsymbol{x}_{0} \\
\boldsymbol{0} & 1
\end{pmatrix}
$, where $\boldsymbol{x}_{0} \in GF(q_{n}^{2})^{d_{n}}$ is any nonzero
vector. Let $\phi \in U_{n}$ be an arbitrary element. Then 
$\rho ( \phi ) =
\begin{pmatrix}
B & \boldsymbol{y} \\
\boldsymbol{0} & 1
\end{pmatrix}
$ for some $\boldsymbol{y} \in GF(q_{n}^{2})^{d_{n}}$ and some strictly
upper triangular matrix $B \in GL(d_{n}, q_{n}^{2})$. Clearly there exist
$A_{1}$, $A_{2} \in SL(d_{n},q_{n}^{2})$ such that
\[
\begin{pmatrix}
\boldsymbol{1} & - \boldsymbol{y} \\
\boldsymbol{0} & 1
\end{pmatrix}
=
\begin{pmatrix}
\boldsymbol{1} & A_{1} \boldsymbol{x}_{0} + A_{2} \boldsymbol{x}_{0} \\
\boldsymbol{0} & 1
\end{pmatrix}
=
\begin{pmatrix}
\boldsymbol{1} & A_{1} \boldsymbol{x}_{0} \\
\boldsymbol{0} & 1
\end{pmatrix}
\begin{pmatrix}
\boldsymbol{1} & A_{2} \boldsymbol{x}_{0} \\
\boldsymbol{0} & 1
\end{pmatrix}
.
\]
Hence there exist $\psi_{1}$, $\psi_{2} \in SU(2d_{n},q_{n}) \cap
\Gamma_{n}$ such that
$\rho ( \psi_{1} \theta \psi_{1}^{-1} \psi_{2} \theta \psi_{2}^{-1}
\phi ) =
\begin{pmatrix}
B & \boldsymbol{0} \\
\boldsymbol{0} & 1
\end{pmatrix}
$ and hence $\psi_{1} \theta \psi_{1}^{-1} \psi_{2} \theta \psi_{2}^{-1}
\phi \in SU(2d_{n},q_{n})$. Thus we can ``uniformly generate'' $U_{n}$
using the element $\theta$. It follows that there exists an element
$g \in \prod_{n}SU(2d_{n}+1,q_{n})$ such that $\prod_{n} U_{n} \leqslant
\langle \prod_{n} SU(2d_{n},q_{n}), g \rangle$.

Now let $\psi \in N_{n} \cap SU(2d_{n},q_{n})$ correspond to the
permutation $(\, e_{1} \, f_{1} \,) \dots (\, e_{d_{n}} \, f_{d_{n}} \,)$.
Then $V_{n} = \psi U_{n} \psi^{-1}$ is the unipotent subgroup of strictly
lower triangular matrices of $SU(2d_{n}+1,q_{n})$; and we also have that
$\prod_{n} V_{n} \leqslant \langle \prod_{n} SU(2d_{n},q_{n}),g
\rangle$.

We can regard $SU(3,q_{n})$ as the subgroup of $SU(2d_{n}+1,q_{n})$
consisting of those elements $\pi$ such that $\pi (e_{i}) = e_{i}$
and $\pi(f_{i}) = f_{i}$ for all $1 \leq i \leq d_{n}-1$. Now let
$h \in \prod_{n} H_{n}$ be an arbitrary element. Then there exists
$g \in \prod_{n} \left( H_{n} \cap SU(2d_{n}, q_{n}) \right)$ such that
$hg \in \prod_{n} \left(H_{n} \cap SU(3,q_{n}) \right)$.
Consequently, in
order to show that $\prod_{n} H_{n}$ is contained in $\langle
\prod_{n} SU(2d_{n},q_{n}),g \rangle$, it is enough to
show that $\prod \left( H_{n} \cap SU(3,q_{n}) \right) \leqslant
\langle \prod_{n} SU(2d_{n},q_{n}) , g \rangle$. To accomplish this,
we shall use a slightly modified form of \cite[pp. 239--242]{ca1}.
For the rest of this proof, we shall write the elements of
$SU(3,q_{n})$ as $3 \times 3$-matrices with respect to the ordered
basis $( e_{d_{n}}, w, f_{d_{n}} )$. Fix an element
$\epsilon \in GF(q_{n}^{2})$ such that $\epsilon \bar{\epsilon} = -1$.
Suppose that $\lambda$, $t \in GF(q_{n}^{2})$ satisfy
$\lambda^{-1} + \bar{\lambda}^{-1} = t \bar{t}$. Then the matrices
\[
A_{1} =
\begin{pmatrix}
1 & \epsilon^{-1} \lambda t & \lambda \\
0 & 1 & \epsilon \bar{\lambda} \bar{t} \\
0 & 0 & 1
\end{pmatrix}
\text{ and } A_{2} =
\begin{pmatrix}
1 & \epsilon^{-1} \bar{\lambda} t & \lambda \\
0 & 1 & \epsilon \lambda \bar{t} \\
0 & 0 & 1
\end{pmatrix}
\]
are elements of $U_{n} \cap SU(3,q_{n})$, and the matrix
\[
B =
\begin{pmatrix}
1 & 0 & 0 \\
- \epsilon \bar{t} & 1 & 0 \\
\bar{\lambda}^{-1} & - \epsilon^{-1} t & 1
\end{pmatrix}
\]
is an element of $V_{n} \cap SU(3,q_{n})$. The product of these matrices
is
\[
A_{1} B A_{2} =
\begin{pmatrix}
0 & 0 & \lambda \\
0 & - \lambda^{-1} \bar{\lambda} & 0 \\
\bar{\lambda}^{-1} & 0 & 0
\end{pmatrix}
.
\]
Let $L_{n}$ be the subset of $GF(q_{n}^{2})^{*}$ consisting of those
elements $\lambda$ such that there exists $t \in GF(q_{n}^{2})$
such that $\lambda^{-1} + \bar{\lambda}^{-1} = t \bar{t}$. By
\cite[13.7.3]{ca1}, each $\lambda \in GF(q_{n}^{2})^{*}$ can be
expressed as $\lambda = \lambda_{1} \bar{\lambda}_{2}^{-1}$ for
some $\lambda_{1}$, $\lambda_{2} \in L_{n}$. Hence we can ``uniformily
generate'' each element of $H_{n} \cap SU(3,q_{n})$ via the equation
\[
\begin{pmatrix}
\lambda & 0 & 0 \\
0 & \lambda^{-1} \bar{\lambda} & 0 \\
0 & 0 & \bar{\lambda}^{-1}
\end{pmatrix}
=
\begin{pmatrix}
0 & 0 & \lambda_{1} \\
0 & - \lambda_{1}^{-1} \bar{\lambda}_{1} & 0 \\
\bar{\lambda}_{1}^{-1} & 0 & 0
\end{pmatrix}
\begin{pmatrix}
0 & 0 & \lambda_{2} \\
0 & - \lambda_{2}^{-1} \bar{\lambda}_{2} & 0 \\
\bar{\lambda}_{2}^{-1} & 0 & 0
\end{pmatrix}
.
\]
\end{proof}

We can now easily obtain the following result.

\begin{Cor} \label{C:unit}
Suppose that $\langle S_{n} \mid n \in \mathbb{N} \rangle$ is a 
sequence of finite simple unitary groups such that there does {\em not\/}
exist an infinite subset $I$ of $\mathbb{N}$ for which conditions
\ref{T:count}(1) and \ref{T:count}(2) are satisfied. Then
$c( \prod_{n} S_{n}) > \omega$.
\end{Cor}
\begin{flushright}
$\square$
\end{flushright} 

\subsection{Orthogonal groups} \label{SS:orth}
In this subsection, we shall consider products of finite orthogonal
groups. First consider the case when each group has the form
$\Omega^{+}(2d,q)$. Then the corresponding orthogonal space has
a normal basis $\boldsymbol{e} \sphat \, \boldsymbol{f}$. Arguing
as in Subsection \ref{SS:symp}, we obtain the following result.

\begin{Thm} \label{T:orth1}
Suppose that $\langle \Omega^{+}(2d_{n},q_{n}) \mid n \in \mathbb{N}
\rangle$ is a sequence of orthogonal groups which satisfies the
following conditions.
\begin{enumerate}
\item[(1)] $d_{n} \geq 3$ for each $n \in \mathbb{N}$.
\item[(2)] There does {\em not\/} exist an infinite subset $I$ of
$\mathbb{N}$ and an integer $d$ such that
\begin{enumerate} 
\item $d_{n} = d$ for all $n \in I$; and
\item if $n$, $m \in I$ and $n < m$, then $q_{n} < q_{m}$.
\end{enumerate}
\end{enumerate}
Then $c(\prod_{n} \Omega^{+}(2d_{n},q_{n})) > \omega$.
\end{Thm}
\begin{flushright}
$\square$
\end{flushright}

Now we shall consider products of the form
$\prod_{n} \Omega(2d_{n}+1,q_{n})$, where $d_{n} \geq 2$ for each 
$n \in \mathbb{N}$. Fix some $n \in \mathbb{N}$. Let $Q$ be the quadratic
form on the corresponding orthogonal space, and let
$(u,v) = Q(u+v) - Q(u) - Q(v)$ be the associated bilinear map. We
can suppose that there exists a basis
\[
\boldsymbol{e} \sphat \, \boldsymbol{f} \sphat \, (w) =
(e_{i} \mid 1 \leq i \leq d_{n} ) \sphat \,
(f_{i} \mid 1 \leq i \leq d_{n}) \sphat \, (w)
\]
of the orthogonal space such that
\[
(e_{i}, f_{j}) = \delta_{ij} \text{ and } (e_{i}, e_{j}) =
(f_{i}, f_{j}) = Q(e_{i}) = Q(f_{i}) = 0
\]
for all $1 \leq i,j \leq d_{n}$; and
\[
Q(w) = 1 \text{ and } (w,e_{i}) = (w,f_{i}) = 0
\]
for all $1 \leq i \leq d_{n}$. Then we can regard $\Omega^{+}(2d_{n},q_{n})$
as the subgroup of $\Omega(2d_{n}+1,q_{n})$ consisting of the elements
$\pi$ such that $\pi$ stabilises the subspace
$\langle e_{i}, f_{i} \mid 1 \leq i \leq d_{n} \rangle$ setwise and
$\pi(w) = w$. Clearly this situation is very similar to
that which we considered in Subsection \ref{SS:unit}. The main difference
is that the Weyl group gets larger in the passage from
$\Omega^{+}(2d_{n},q_{n})$ to $\Omega(2d_{n}+1,q_{n})$. The Weyl group
$W_{n}$ of $\Omega(2d_{n}+1,q_{n})$ is
$\mathbb{Z}_{2}^{d_{n}} \rtimes Sym(d_{n})$, acting on the set \\
$\{ e_{i},f_{i} \mid 1 \leq i \leq d_{n} \}$ with blocks of imprimitivity
$\{ e_{i},f_{i} \mid 1 \leq i \leq d_{n} \}$. The Weyl group of
$\Omega^{+}(2d_{n},q_{n})$ is the subgroup $W_{n}^{+}$ of $W_{n}$
consisting of the even permutations of
$\{ e_{i},f_{i} \mid 1 \leq i \leq d_{n} \}$. But this point has
already been dealt with during our treatment of the symplectic
groups in Subsection \ref{SS:symp}. Hence we can easily obtain the
following result.

\begin{Thm} \label{T:orth2}
Suppose that $d_{n} \geq 2$ for all $n \in \mathbb{N}$. Then
$\prod_{n} \Omega(2d_{n}+1,q_{n})$ is finitely generated over
$\prod_{n} \Omega^{+}(2d_{n},q_{n})$.
\end{Thm}
\begin{flushright}
$\square$
\end{flushright}

Finally we shall consider products of the form
$\prod_{n} \Omega^{-}(2d_{n}+2,q_{n})$, where $d_{n} \geq 3$ for 
each $n \in \mathbb{N}$. In this case, there exists a basis
$\boldsymbol{e} \sphat \, \boldsymbol{f} \sphat \, (w,z)$ of
the corresponding orthogonal space such that
\[
(e_{i},f_{j}) = \delta_{ij} \text{ and }
(e_{i},e_{j}) = (f_{i},f_{j})= Q(e_{i}) = Q(f_{i}) =0
\]
for all $1 \leq i,j \leq d_{n}$; and
\[
(w,e_{i}) = (w,f_{i}) = (z,e_{i})= (z,f_{i})=0
\]
for all $1 \leq i \leq d_{n}$; and the subspace $\langle w,z \rangle$
does not contain any singular vectors. So we can regard
$\Omega^{+}(2d_{n},q_{n})$ as the subgroup of
$\Omega^{-}(2d_{n}+2,q_{n})$ consisting of the elements $\pi$
such that $\pi$ stabilises the subspace
$\langle e_{i}, f_{i} \mid 1 \leq i \leq d_{n} \rangle$ setwise,
$\pi(w) = w$ and $\pi(z) = z$.

\begin{Thm} \label{T:orth3}
Suppose that $d_{n} \geq 3$ for each $n \in \mathbb{N}$. Then
$\prod_{n} \Omega^{-}(2d_{n}+2,q_{n})$ is finitely generated
over $\prod_{n} \Omega^{+}(2d_{n},q_{n})$.
\end{Thm}

\begin{proof}
As before, we shall make use of the Bruhat decomposition
\[
\Omega^{-}(2d+2,q) = \underset{w \in W}{\bigcup} BwB ,
\]
where $B$ is a Borel subgroup and $W$ is the Weyl group. Fix some
$n \in \mathbb{N}$. For each $1 \leq i \leq d_{n}$, let
$E_{i} = \langle e_{1}, \dots ,e_{i} \rangle$. Then the
stabiliser $B_{n}$ of the flag of totally singular subspaces
\[
E_{1} \leqslant E_{2} \leqslant \dots \leqslant E_{d_{n}}
\]
is a Borel subgroup of $\Omega^{-}(2d_{n}+2,q_{n})$. Let $N_{n}$ be the
subgroup of $\Omega^{-}(2d_{n}+2,q_{n})$ which stabilises the polar
frame $\{ \langle e_{i} \rangle , \langle f_{i} \rangle \mid
1 \leq i \leq d_{n} \}$. Then the Weyl group of
$\Omega^{-}(2d_{n}+2,q_{n})$ is $W_{n} = N_{n}/B_{n} \cap N_{n}$. Once
again, $W_{n}$ is $\mathbb{Z}_{2}^{d_{n}} \rtimes Sym(d_{n})$ acting on
the set $\{ e_{i},f_{i} \mid 1 \leq i \leq d_{n} \}$ with blocks of
imprimitivity $\{ e_{i},f_{i} \}$ for $1 \leq i \leq d_{n}$. As before,
the main point is to show that there exists a subgroup $G$ of
$\prod_{n} \Omega^{-}(2d_{n}+2,q_{n})$ such that
\begin{enumerate}
\item $G$ is finitely generated over $\prod_{n} \Omega^{+}(2d_{n},q_{n})$,
and
\item $\prod_{n} B_{n} \leqslant G$.
\end{enumerate}

Let $U_{n}$ be the subgroup of unipotent elements of $B_{n}$ and let
$H_{n} = B_{n} \cap N_{n}$; so that $B_{n} = U_{n} \rtimes H_{n}$.
First we shall show that there exists an element
$g_{0} \in \prod_{n} \Omega^{-}(2d_{n}+2,q_{n})$ such that
$\prod_{n} U_{n} \leqslant G_{0} = \langle
\prod_{n} \Omega^{+}(2d_{n},q_{n}), g_{0} \rangle$. Note that if
$\pi \in U_{n}$, then there exist vectors
$\boldsymbol{x}$, $\boldsymbol{y} \in E_{d_{n}}$ such that
$\pi(w) = w + \boldsymbol{x}$ and
$\pi(z) = z+ \boldsymbol{y}$. Let
$E_{d_{n}}^{+} = \langle E_{d_{n}}, w,z \rangle$ and let 
$\Gamma_{n}$ be the setwise stabiliser of $E_{d_{n}}^{+}$ in
$\Omega^{-}(2d_{n}+2,q_{n})$. Let $\rho : \Gamma_{n} \to
GL(E_{d_{n}}^{+})$ be the restriction map. We regard $GL(E_{d_{n}}^{+})$
as a group of matrices with respect to the ordered basis
$(e_{1}, \dots ,e_{d_{n}},w,z)$. So for each $A \in SL(d_{n},q_{n})$,
we have that 
\[
\begin{pmatrix}
A & \boldsymbol{0} \\
\boldsymbol{0} & \boldsymbol{1}
\end{pmatrix}
\in \rho \left[ \Omega^{+}(2d_{n},q_{n}) \cap \Gamma_{n} \right] .
\]
Arguing as in the proof of Theorem \ref{T:unit2}, we see that the
existence of a suitable element $g_{0} \in
\prod_{n} \Omega^{-}(2d_{n}+2,q_{n})$ is a consequence of the
following easy observation.

\begin{Claim}
Suppose that $d \geq 3$. Let $S = SL(d,q)$ and $V = V(d,q)$.
Let $S$ act on $V \times V$ via the action
$A( \boldsymbol{x}, \boldsymbol{y} ) =
(A \boldsymbol{x}, A \boldsymbol{y} )$. Suppose that
$\boldsymbol{a}$, $\boldsymbol{b} \in V$ are linearly independent.
Then for all $(\boldsymbol{x}, \boldsymbol{y}) \in V \times V$,
there exist $A$, $B \in S$ such that 
$(\boldsymbol{x}, \boldsymbol{y}) = 
A(\boldsymbol{a}, \boldsymbol{b}) +
B(\boldsymbol{a}, \boldsymbol{b})$.
\end{Claim}
\begin{flushright}
$\square$
\end{flushright}

Finally we shall show that there exists an element 
$g_{1} \in \prod_{n} \Omega^{-}(2d_{n}+2,q_{n})$ such that
$\prod_{n} H_{n} \leqslant G_{1} = \langle G_{0},g_{1} \rangle$.
We shall regard $\Omega^{-}(4,q_{n})$ as the subgroup of 
$\Omega^{-}(2d_{n}+2,q_{n})$ consisting of the elements $\pi$
such that $\pi$ stabilises the subspace
$\langle e_{d_{n}}, f_{d_{n}}, w, z \rangle$ setwise and such that
$\pi(e_{i}) = e_{i}$ and $\pi(f_{i}) = f_{i}$ for all
$1 \leq i \leq d_{n}-1$. Since
$\prod_{n} \Omega^{+}(2d_{n},q_{n}) \leqslant G_{0}$, it is enough
to find an element $g_{1}$ such that
$\prod_{n} \left( H_{n} \cap \Omega^{-}(4,q_{n}) \right)
\leqslant \langle G_{0}, g_{1} \rangle$. We shall make use of
the fact that
\[
\Omega^{-}(4,q_{n}) \simeq SL(2,q_{n}^{2})/ \{ \pm \boldsymbol{1} \}.
\]
(For example, see \cite[12.42]{ta}.) Let $p_{n} = 
\text{char} (GF(q_{n}))$. Then $U_{n} \cap \Omega^{-}(4,q_{n})$ is a
group of order $q_{n}^{2}$, and hence is a Sylow $p_{n}$-subgroup
of $\Omega^{-}(4,q_{n})$. It is easily checked that $SL(2,q_{n}^{2})$
is ``uniformly generated'' by the two subgroups
$UT(2,q_{n}^{2})$ and $LT(2,q_{n}^{2})$, consisting of the strictly upper
triangular matrices and strictly lower triangular matrices of
$SL(2,q_{n}^{2})$. (See \cite[6.4.4]{ca1}.) Since $UT(2,q_{n}^{2})$ and
$LT(2,q_{n}^{2})$ are Sylow $p_{n}$-subgroups of $SL(2,q_{n}^{2})$,
it follows that there exists an element
$g_{1} \in \prod_{n} \Omega^{-}(2d_{n}+2,q_{n})$ such that
$\prod_{n} \Omega^{-}(4,q_{n}) \leqslant 
\langle G_{0}, g_{1} \rangle$.   
\end{proof}

\begin{Cor} \label{C:orth}
Suppose that $\langle S_{n} \mid n \in \mathbb{N} \rangle$ is a sequence
of finite simple orthogonal groups such that there does {\em not\/}
exist an infinite subset $I$ of $\mathbb{N}$ for which conditions
\ref{T:count}(1) and \ref{T:count}(2) are satisfied. Then
$c( \prod_{n} S_{n} ) > \omega$.
\end{Cor} 

\subsection{Conclusion} \label{SS:conc}
We can now complete the proof of Theorem \ref{T:class}. Suppose that
$\langle S_{n} \mid n \in \mathbb{N} \rangle$ is sequence of finite
simple nonabelian groups such that there does {\em not\/} exist
an infinite subset $I$ of $\mathbb{N}$ for which conditions
\ref{T:count}(1) and \ref{T:count}(2) are satisfied. Let
$G = \prod_{n} S_{n}$. Let
\begin{itemize}
\item $\mathcal{C}_{0}$ be the set of 26 sporadic finite simple groups,
\item $\mathcal{C}_{1}$ be the set of finite simple alternating groups,
\item $\mathcal{C}_{2}$ be the set of finite simple projective special
linear groups,
\item $\mathcal{C}_{3}$ be the set of finite simple symplectic groups,
\item $\mathcal{C}_{4}$ be the set of finite simple unitary groups,
\item $\mathcal{C}_{5}$ be the set of finite simple orthogonal groups,
and
\item $\mathcal{C}_{6}$ be the set of finite simple groups of Lie types
$E_{6}$, $E_{7}$, $E_{8}$, $F_{4}$, $G_{2}$, $^{2}E_{6}$,
$^{2}B_{2}$, $^{2}G_{2}$, $^{2}F_{4}$ and $^{3}D_{4}$.
\end{itemize}
By the classification of the finite simple groups, each finite simple
nonabelian group lies in one of the above sets. Some groups lie in more
than one of these sets. For example, $Alt(8) \simeq PSL(4,2)$. For the rest
of this argument, we shall suppose that we have slightly modified the
above sets so that they yield a partition of the finite simple nonabelian
groups.

For each $0 \leq i \leq 6$, let
$J_{i} = \{ n\in \mathbb{N} \mid S_{n} \in \mathcal{C}_{i} \}$ and let
$P_{i} = \prod_{n \in J_{i}} S_{n}$. Then 
$G = \prod_{i = 0}^{6} P_{i}$. Using Proposition \ref{P:gen}, it is 
enough to show that for each $0 \leq i \leq 6$, either
$c(P_{i}) > \omega$ or $P_{i}$ is finite. If $1 \leq i \leq 5$,
this has been proved in Sections \ref{S:alt}, \ref{S:lin} and
\ref{S:class}.
And if $i = 0$, this is an immediate consequence of Proposition
\ref{P:finite}. Finally consider
$P_{6}$. Our hypothesis on 
$\langle S_{n} \mid n \in \mathbb{N} \rangle$ implies that there
exists a finite set of simple groups
$\mathcal{F} \subseteq \mathcal{C}_{6}$ such that $S_{n} \in \mathcal{F}$
for all $n \in J_{6}$. So the result once again follows from
Proposition \ref{P:finite}. This completes
the proof of Theorem \ref{T:class}.  

\section{A consistency result} \label{S:con}
In this section, we shall prove Theorem \ref{T:con}. Our notation follows
that of Kunen \cite{ku}. Thus if $\mathbb{P}$ is a notion of forcing and
$p$, $q \in \mathbb{P}$, then $q \leq p$ means that $q$ is a strengthening
of $p$. If $V$ is the ground model, then we denote the generic extension
by $V^{\mathbb{P}}$ when we do not want to specify a particular generic
filter $G \subseteq \mathbb{P}$.

\begin{Def} \label{D:laver}
A notion of forcing $\mathbb{P}$ is said to have the {\em Laver property\/}
if the following holds. Suppose that
\begin{enumerate}
\item $\langle A_{n} \mid n \in \mathbb{N} \rangle$ is a sequence of
finite sets;
\item $f: \mathbb{N} \to \mathbb{N}$ is a function such that $f(n) \geq 1$ for all
$n \in \mathbb{N}$ and
$f(n) \to
\infty$ as $n \to \infty$;
\item $p \in \mathbb{P}$, $\tilde{g}$ is a $\mathbb{P}$-name and
$p \Vdash \tilde{g} \in \prod_{n} A_{n}$.
\end{enumerate}
Then there exists $q \leq p$ and a sequence $\langle B_{n} \mid n \in \mathbb{N} 
\rangle$ such that
\begin{enumerate}
\item[(a)] $B_{n} \subseteq A_{n}$ and $\left|B_{n} \right| \leq f(n)$;
\item[(b)] $q \Vdash \tilde{g} \in \prod_{n} B_{n}$.
\end{enumerate}
\end{Def}

Theorem \ref{T:con} is an immediate consequence of the following two
results.

\begin{Thm} \label{T:cover}
Suppose that $V \vDash CH$, and that
$\langle \mathbb{P}_{\alpha}, \tilde{\mathbb{Q}}_{\alpha} \mid \alpha < \omega_{2}
\rangle$
is a countable support iteration of proper notions of forcing such that
for all $\alpha < \omega_{2}$
\begin{enumerate}
\item[(1)] $\Vdash_{\alpha} \tilde{\mathbb{Q}}_{\alpha}$ has the cardinality
of the continuum; and
\item[(2)] $\Vdash_{\alpha} \tilde{\mathbb{Q}}_{\alpha}$ has the Laver property.
\end{enumerate}
Then in $V^{\mathbb{P}_{\omega_{2}}}$, $c \left( \prod_{n} G_{n} \right)
\leq \omega_{1}$ for {\em every\/} sequence 
$\langle G_{n} \mid n \in \mathbb{N} \rangle$ of nontrivial finite 
groups.
\end{Thm}

\begin{Thm} \label{T:sym}
Suppose that $V \vDash CH$. Then there exists a countable support
iteration
$\langle \mathbb{P}_{\alpha}, \tilde{\mathbb{Q}}_{\alpha} \mid
\alpha < \omega_{2} \rangle$
of proper notions of forcing such that
\begin{enumerate}
\item[(a)] $\langle \mathbb{P}_{\alpha}, \tilde{\mathbb{Q}}_{\alpha} \mid
\alpha < \omega_{2} \rangle$ satisfies conditions \ref{T:cover}(1)
and \ref{T:cover}(2); and
\item[(b)] $V^{\mathbb{P}_{\omega_{2}}} \vDash
c(Sym(\mathbb{N})) = \omega_{2} = 2^{\omega}$.
\end{enumerate}
\end{Thm}

First we shall prove Theorem \ref{T:cover}.

\begin{Def} \label{D:cover}
Let $\langle G_{n} \mid n \in \mathbb{N} \rangle$ be a sequence of
nontrivial finite groups.
\begin{enumerate}
\item A {\em cover\/} is a function 
$c: \mathbb{N} \to \left[ \underset{n \in \mathbb{N}}{\bigcup}G_{n} \right]^{< \omega}$
such that for all $n \in \mathbb{N}$
\begin{enumerate}
\item $\emptyset \ne c(n) \subseteq G_{n}$;
\item the identity element $1_{G_{n}} \in c(n)$;
\item if $a \in c(n)$, then $a^{-1} \in c(n)$.
\end{enumerate}
\item If $g = \langle g(n) \rangle_{n} \in \prod_{n} G_{n}$, then
$c$ {\em covers\/} $g$ if $g(n) \in c(n)$ for all $n \in \mathbb{N}$.
\item If $c$ is a cover and $f: \mathbb{N} \to \mathbb{N}$, then $c$ is
an $f$-{\em cover\/} if $|c(n)| \leq f(n)$ for all $n \in \mathbb{N}$.
\item If $c_{1}$ and $c_{2}$ are covers, then the cover $c_{1} * c_{2}$
is defined by
\[
(c_{1}*c_{2})(n) = \{ ab , (ab)^{-1} \mid a \in c_{1}(n), b \in c_{2}(n) \}.
\]
\end{enumerate}
\end{Def}

\begin{Lem} \label{L:cover}
If $c_{1}$ is an $f_{1}$-cover and $c_{2}$ is an $f_{2}$-cover,
then $c_{1}*c_{2}$ is a $2f_{1} f_{2}$-cover.
\end{Lem}

\begin{proof}
Obvious.
\end{proof}

It is perhaps worth mentioning that $*$ is generally {\em not\/} an
associative operation on the set of covers of
$\prod_{n} G_{n}$.

\begin{Def} \label{D:close}
If $C$ is a set of covers of $\prod_{n} G_{n}$, then its {\em closure\/}
$c\ell (C)$ is the least set of covers satisfying
\begin{enumerate}
\item $C \subseteq c\ell (C)$; and
\item if $d_{1}$, $d_{2} \in c\ell (C)$ then $d_{1}*d_{2} \in c\ell (C)$.
\end{enumerate}
\end{Def}

\begin{Lem} \label{L:close}
Suppose that $C$ is a set of covers of $\prod_{n} G_{n}$. Then
\[
\{ g \in \prod_{n} G_{n} \mid
\text{There exists } d \in c\ell (C) \text{ such that } g \text{ is
covered by } d \}
\]
is a subgroup of $\prod_{n} G_{n}$.
\end{Lem}

\begin{proof}
Easy.
\end{proof}

 From now on, let $f: \mathbb{N} \to \mathbb{N}$ be the function defined
by $f(n) = 2^{n+2}$ for all $n \in \mathbb{N}$.

\begin{Lem} \label{L:proper}
Suppose that $\langle G_{n} \mid n \in \mathbb{N} \rangle$ is a sequence
of finite groups such that $\left| G_{n} \right| \geq 2^{(n+2)^{2}}$
for all $n \in \mathbb{N}$. If $C$ is a countable set of $f$-covers of
$\prod_{n} G_{n}$, then
\[
\{ g \in \prod_{n} G_{n} \mid \text{There exists } d \in c\ell (C)
\text{ such that } g \text{ is covered by } d \}
\]
is a {\em proper\/} subgroup of $\prod_{n} G_{n}$.
\end{Lem}

\begin{proof}
Suppose that $d \in c\ell (C)$ is an $m$-fold $*$-product of
$c_{1}, \dots , c_{m} \in C$ in some order. (Remember that $*$
is not an associative operation.) Then Lemma \ref{L:cover}
implies that $d$ is a $2^{m-1}f^{m}$-cover. So we can 
enumerate $c\ell (C) = \{ d_{n} \mid n \in \mathbb{N} \}$ in
such a way that $d_{n}$ is a $\phi_{n}$-cover for all $n \in \mathbb{N}$, where
$\phi_{n} = 2^{n}f^{n+1}$. In particular,
\[
\left| d_{n}(n) \right| \leq 2^{n}f(n)^{n+1} = 2^{n^{2}+4n+2} <
\left| G_{n} \right|.
\]
Hence there exists $g = \langle g(n) \rangle_{n} \in \prod_{n} G_{n}$
such that $g(n) \in G_{n} \smallsetminus d_{n}(n)$ for all
$n \in \mathbb{N}$. Clearly $g$ is not covered by any element
$d \in c\ell (C)$.
\end{proof}

\begin{proof}[Proof of Theorem \ref{T:cover}]
Suppose that $V \vDash CH$ and that
$\langle \mathbb{P}_{\alpha} , \tilde{\mathbb{Q}}_{\alpha}
\mid \alpha < \omega_{2} \rangle$ is a countable support iteration
of proper notions of forcing such that for all $\alpha < \omega_{2}$
\begin{enumerate}
\item $\Vdash_{\alpha} \tilde{\mathbb{Q}}_{\alpha}$ has the cardinality
of the continuum; and
\item $\Vdash_{\alpha} \tilde{\mathbb{Q}}_{\alpha}$ has the Laver property.
\end{enumerate}

 From now on, we shall work inside $V^{\mathbb{P}_{\omega_{2}}}$. Let
$\langle G_{n} \mid n \in \mathbb{N} \rangle$ be a sequence of nontrivial
finite groups. First suppose that there exists an infinite subset $I$
of $\mathbb{N}$ and a finite group $G$ such that $G_{n} = G$ for all
$n \in I$. By Lemma \ref{L:normal} and Theorems
\ref{T:tits} and \ref{T:exact},
$c \left( \prod_{n} G_{n} \right) \leq c \left( \prod_{n \in I} G_{n} \right)
\leq \omega_{1}$.
Hence we can assume that no such subset $I$ of $\mathbb{N}$ exists. Then
there exists an infinite subset $J = \{ j_{n} \mid n \in \mathbb{N} \}$
of $\mathbb{N}$ such that $\left| G_{j_{n}} \right| \geq 2^{(n+2)^{2}}$
for all $n \in \mathbb{N}$. By Lemma \ref{L:normal},
$c \left( \prod_{n} G_{n} \right) \leq c \left( \prod_{n \in J} G_{n} \right)$.
To simplify notation, we shall suppose that 
$\left| G_{n} \right| \geq 2^{(n+2)^{2}}$ for all $n \in \mathbb{N}$.

Since the sequence $\langle G_{n} \mid n \in \mathbb{N} \rangle$ of
finite groups can be coded by a real number,
there exists $\alpha < \omega_{2}$ such that 
$\langle G_{n} \mid n \in \mathbb{N} \rangle \in V^{\mathbb{P}_{\alpha}}$.
By Shelah III 4.1 \cite{sh-b}, $V^{\mathbb{P}_{\alpha}} \vDash CH$. Let
$\{ c_{\beta} \mid \beta < \omega_{1} \}$ be an enumeration of the
$f$-covers $c \in V^{\mathbb{P}_{\alpha}}$ of $\prod_{n} G_{n}$. For
each $\gamma < \omega_{1}$, let
$C_{\gamma} = \{ c_{\beta} \mid \beta < \gamma \}$ and define
\[
H_{\gamma} = \{ g \in \prod_{n} G_{n} \mid \text{There exists } 
d \in c\ell (C_{\gamma} ) \text{ such that } g \text{ is covered by }
d \}.
\]
By Lemma \ref{L:proper}, $H_{\gamma}$ is a proper subgroup of 
$\prod_{n} G_{n}$ for all $\gamma < \omega_{1}$. Thus it suffices to
show that $\prod_{n} G_{n} = \underset{\gamma < \omega_{1}}{\bigcup}
H_{\gamma}$.

Let $g \in \prod_{n} G_{n}$ be any element. By Shelah \cite[VI Section 2]{sh-b} 
and \cite[Appendix]{sh}, the Laver property is preserved by countable support
iterations of proper notions of forcing. This implies that there
exists a sequence 
$\langle B_{n} \mid n \in \mathbb{N} \rangle \in V^{\mathbb{P}_{\alpha}}$
such that
\begin{enumerate}
\item[(a)] $B_{n} \subseteq G_{n}$ and $\left| B_{n} \right|
\leq 2^{n}$; and
\item[(b)] $g(n) \in B_{n}$ for all $n \in \mathbb{N}$.
\end{enumerate}
Define the function $c$ by
\[
c(n) = B_{n} \cup \{ a^{-1} \mid a \in B_{n} \} \cup \{ 1_{G_{n}} \}
\]
for all $n \in \mathbb{N}$. Then $c \in V^{\mathbb{P}_{\alpha}}$ is an
$f$-cover of $\prod_{n} G_{n}$, and so $c = c_{\beta}$ for
some $\beta < \omega_{1}$. Hence 
$g \in \underset{\gamma < \omega_{1}}{\bigcup}H_{\gamma}$. 
\end{proof}

The rest of this section will be devoted to the proof of Theorem
\ref{T:sym}. Each of the notions of forcing which we shall use
in our iteration will satisfy Axiom A. It is well-known that if
$\mathbb{P}$ satisfies Axiom A, then $\mathbb{P}$ is proper.
(For example, see \cite[p.101]{j}. )

\begin{Def} \label{D:axiom}
A notion of forcing $\mathbb{P}$ satisfies {\em Axiom A\/} if there
is a collection $\{ \leq_{n} \mid n \in \omega \}$ of partial
orderings of $\mathbb{P}$ which satisfies the following conditions.
\begin{enumerate}
\item[(1)] $p \leq_{0} q$ if{f} $p \leq q$.
\item[(2)] If $p \leq_{n+1} q$, then $p \leq_{n} q$.
\item[(3)] If $\langle p_{n} \mid n \in \omega \rangle$ is a sequence
such that $p_{n+1} \leq_{n} p_{n}$ for all $n \in \omega$, then
there exists $q \in \mathbb{P}$ such that $q \leq_{n} p_{n}$ for
all $n \in \omega$.
\item[(4)] For each $p \in \mathbb{P}$, $n \in \omega$ and ordinal
name $\tilde{\alpha}$, there exists a countable set $B$ and a 
condition $q \in \mathbb{P}$ such that $q \leq_{n} p$ and
$q \Vdash \tilde{\alpha} \in B$.
\end{enumerate}
\end{Def}

\begin{Def} \label{D:poset}
Fix a partition 
$\{ I_{n} \mid n \in \mathbb{N} \}$ of $\mathbb{N}$ into infinitely many
finite subsets such that the following conditions hold.
\begin{enumerate}
\item $\left| I_{n} \right| \geq 2$ for all $n \in \mathbb{N}$.
\item For each $t \geq 2$, there exist infinitely many $n \in \mathbb{N}$
such that $\left| I_{n} \right| = t$.
\item If $n < m$, then $\max(I_{n}) < \min(I_{m})$. (Thus each $I_{n}$
consists of a finite set of consecutive integers.)
\end{enumerate}
The notion of forcing $\mathbb{B}$ consists of all functions $p$ such that
\begin{enumerate}
\item[(a)] there exists a subset $J$ of $\mathbb{N}$ such that 
$\dom p = \underset{n \in J}{\bigcup}I_{n}$;
\item[(b)]if $n \in J$, then $p \res I_{n} \in Sym(I_{n})$;
\item[(c)] if $t \geq 2$, then there exist infinitely many 
$n \in \mathbb{N} \smallsetminus J$ such that 
$\left| I_{n} \right| = t$.
\end{enumerate}
If $p$, $q \in \mathbb{B}$, then we define $q \leq p$ if and only if
$q \supseteq p$.
\end{Def}

\begin{Lem} \label{L:poset}
$\mathbb{B}$ satisfies Axiom A and has the Laver property.
\end{Lem}

\begin{proof}
For each $p \in \mathbb{B}$ and $t \geq 2$, let
\[
S^{t}(p) = \{ m \in \mathbb{N} \mid m \notin \dom p \text{ and }
\left| I_{m} \right| = t \} ;
\]
and for each $n \geq 1$, let $S^{t}_{n}(p)$ be the set of the first
$n$ elements of $S^{t}(p)$. If $t \leq 1$ or $n = 0$, let
$S^{t}_{n}(p) = \emptyset$. For each $n \in \omega$, define
a partial ordering $\leq_{n}$ on $\mathbb{B}$ by setting
$q \leq_{n} p$ if and only if
\begin{enumerate}
\item $q \supseteq p$, and
\item for each $t \leq n$, $S^{t}_{n}(p) \subseteq \mathbb{B}
\smallsetminus \dom p$.
\end{enumerate}
Then it is easily checked that the partial orderings
$\{ \leq_{n} \mid n \in \omega \}$ satisfy clauses (1)--(4)
of Definition \ref{D:axiom}. It is also easy to verify that
$\mathbb{B}$ has the Laver property.
\end{proof}

It follows that $\mathbb{B}$ satisfies conditions \ref{T:cover}(1) and \ref{T:cover}(2).
(It is also easily seen that $\mathbb{B}$ is $^{\omega}\omega$-bounding.
However, we shall not need this fact in the proof of Theorem \ref{T:sym}.)
After first introducing some group theoretic notation, we shall explain
the relevance of $\mathbb{B}$ to the problem of computing $c(Sym(\mathbb{N}))$.

\begin{Def} \label{D:perm}
Suppose that $\{ a_{n} \mid n \in \mathbb{N} \}$ is the increasing 
enumeration of the infinite subset $A$ of $\mathbb{N}$. If 
$\pi \in Sym(\mathbb{N})$, then $\pi^{A} \in Sym(A)$ is defined by
$\pi^{A}(a_{n}) = a_{\pi(n)}$ for all $n \in \mathbb{N}$. If
$\Gamma$ is a subgroup of $Sym(\mathbb{N})$, then
$\Gamma^{A} = \{ \pi^{A} \mid \pi \in \Gamma \}$.
\end{Def}

\begin{Def} \label{D:product}
If $g : \mathbb{N} \to \mathbb{N}$ is a strictly increasing function, then
\[
P_{g} = \prod_{n} Sym(g(n) \smallsetminus g(n-1)).
\]
(Here we use the convention that $g(-1) = 0$.)
\end{Def}

\begin{Def} \label{D:growth}
\begin{enumerate} 
\item If $f$, $g : \mathbb{N} \to \mathbb{N}$, then $f \leq^{*} g$ if{f}
there exists $n_{0} \in \mathbb{N}$ such that $f(n) \leq g(n)$ for
all $n \geq n_{0}$.
\item If $g : \mathbb{N} \to \mathbb{N}$ is a strictly increasing function,
then
\[
S_{g} = \langle \pi \in Sym(\mathbb{N} \mid \pi, \pi^{-1} \leq^{*} g
\rangle.
\]
\end{enumerate}
\end{Def}

$\mathbb{B}$ was designed so that the following density condition would
hold.

\begin{Lem} \label{L:density}
Suppose that $g : \mathbb{N} \to \mathbb{N}$ is a strictly increasing function
and that $p \in \mathbb{B}$. Then there exists an infinite subset $A$ of
$\mathbb{N} \smallsetminus \dom p$ such that $p \cup \pi \in \mathbb{B}$
for all $\pi \in P_{g}^{A}$.
\end{Lem}

\begin{proof}
Let $\dom p = \underset{n \in J}{\bigcup}I_{n}$. Then it is easy to
find a suitable set $A$ of the form 
$\underset{n \in K}{\bigcup}I_{n}$, where $K$ is an appropriately
chosen subset of $\mathbb{N} \smallsetminus J$.
\end{proof}

Arguing as in Section 2 of \cite{st2}, we can now easily obtain the 
following result.

\begin{Lem} \label{L:growth}
Let $V \vDash CH$ and let 
$\langle \mathbb{P}_{\alpha}, \tilde{\mathbb{Q}}_{\alpha} \mid
\alpha < \omega_{2} \rangle$ be a countable support iteration
of proper notions of forcing such that for each $\alpha < \omega_{2}$,
$\mathbb{P}_{\alpha} \Vdash | \tilde{\mathbb{Q}}_{\alpha} |
= 2^{\omega}$. 
Suppose that 
$S \subseteq \{ \alpha < \omega_{2} \mid cf(\alpha) = \omega_{1} \}$
is a stationary subset of $\omega_{2}$, and that 
$\tilde{\mathbb{Q}}_{\alpha} = \tilde{\mathbb{B}}$ for all $\alpha \in S$.
(Here $\tilde{\mathbb{B}}$ is the notion of forcing $\mathbb{B}$ in
the generic extension $V^{\mathbb{P}_{\alpha}}$.) Then the following
statements are equivalent in $V^{\mathbb{P}_{\omega_{2}}}$.
\begin{enumerate}
\item[(1)] $c(Sym(\mathbb{N})) = \omega_{1}$.
\item[(2)] It is possible to express $Sym(\mathbb{N}) =
\underset{i < \omega_{1}}{\bigcup}G_{i}$ as the union of a chain
of proper subgroups such that for each strictly increasing function
$g : \mathbb{N} \to \mathbb{N}$, there exists $i < \omega_{1}$ with
$S_{g} \leqslant G_{i}$.
\end{enumerate}
\end{Lem}

\begin{flushright}
$\square$
\end{flushright}

\begin{Def} \label{D:lav}
Laver forcing $\mathbb{L}$ consists of the set of all trees
$T \subseteq {}^{< \omega}\omega$ with the following property.
There exists an integer $k$ such that
\begin{enumerate}
\item if $n < k$, then $\left| T \cap {}^{n} \omega \right| = 1$;
\item if $n \geq k$ and $\eta \in T \cap {}^{n} \omega$, then there
exist infinitely many $i \in \omega$ such that 
$\eta\sphat \, \langle i \rangle \in T$.
\end{enumerate}
If $S$, $T \in \mathbb{L}$, then $S \leq T$ if{f} 
$S \subseteq T$.
\end{Def}

The following result is well-known.

\begin{Lem} \label{L:lav}
\begin{enumerate}
\item[(1)] Suppose that $V \vDash ZFC$. Then there exists a function
$g \in {}^{\mathbb{N}}\mathbb{N} \cap V^{\mathbb{L}}$ such that 
$f \leq^{*} g$ for all $f \in {}^{\mathbb{N}}\mathbb{N} \cap V$.
\item[(2)] $\mathbb{L}$ satisfies Axiom A and has the Laver property.
\end{enumerate}
\end{Lem}

\begin{flushright}
$\square$
\end{flushright}

It is now easy to complete the proof of Theorem \ref{T:sym}.
Let $V \vDash CH$. Define a countable support iteration
$\langle \mathbb{P}_{\alpha}, \tilde{\mathbb{Q}}_{\alpha} \mid
\alpha < \omega_{2} \rangle$ of proper notions of forcing 
with the Laver property inductively as follows. If
$cf(\alpha) = \omega_{1}$, let $\tilde{\mathbb{Q}}_{\alpha}
= \tilde{\mathbb{B}}$. Otherwise, let
$\tilde{\mathbb{Q}}_{\alpha} = \tilde{\mathbb{L}}$. From now on,
we work inside $V^{\mathbb{P}_{\omega_{2}}}$. Clearly
$2^{\omega} = \omega_{2}$.
Suppose that
$c(Sym(\mathbb{N})) = \omega_{1}$. By Lemma \ref{L:growth},
we can express 
$Sym(\mathbb{N}) = \underset{i < \omega_{1}}{\bigcup}G_{i}$
as the union of a chain of proper subgroups such that for
each strictly increasing function $g: \mathbb{N} \to \mathbb{N}$,
there exists $i < \omega_{1}$ with $S_{g} \leqslant G_{i}$.
Lemma \ref{L:lav} implies that there exists a sequence
$\langle g_{\alpha} : \mathbb{N} \to \mathbb{N} \mid \alpha < \omega_{2} 
\rangle$ of strictly increasing functions such that
\begin{enumerate}
\item if $\alpha < \beta < \omega_{2}$, then $g_{\alpha}
\leq^{*} g_{\beta}$; and
\item for all $f: \mathbb{N} \to \mathbb{N}$, there exists $\alpha < \omega_{2}$
such that $f \leq^{*} g_{\alpha}$.
\end{enumerate}
There exists $i < \omega_{1}$ and a cofinal subset $C$ of $\omega_{2}$
such that $S_{g_{\alpha}} \leqslant G_{i}$ for all $\alpha \in C$.
But this means that $G_{i} = Sym(\mathbb{N})$, which is a contradiction.
Hence $c(Sym(\mathbb{N})) = \omega_{2}$.

\section{$Sym(\mathbb{N})$ has property (FA)} \label{S:sym}
In this section, we shall prove that $Sym(\mathbb{N})$ has property (FA).
By Macpherson and Neumann \cite{mn}, $c(Sym(\mathbb{N})) > \omega$.
Also, since every proper normal subgroup of $Sym(\mathbb{N})$ is
countable, $\mathbb{Z}$ is not a homomorphic image of $Sym(\mathbb{N})$.
Thus it is enough to prove the following result.

\begin{Thm} \label{T:amal}
$Sym(\mathbb{N})$ is not a nontrivial free product with amalgamation.
\end{Thm}

Suppose that $Sym(\mathbb{N})$ is a nontrivial free product with
amalgamation. Then there exists a tree $T$ such that
\begin{enumerate}
\item $Sym(\mathbb{N})$ acts without inversion on $T$; and
\item there exists $\pi \in Sym(\mathbb{N})$ such that
$\pi(t) \ne t$ for all $t \in T$.
\end{enumerate}
(See Theorem 7 \cite{se}.) Thus it suffices to prove that whenever
$Sym(\mathbb{N})$ acts without inversion on a tree $T$, then for
every $\pi \in Sym(\mathbb{N})$ there exists a vertex $t \in T$ such
that $\pi(t) = t$. (This also yields a second proof that $\mathbb{Z}$
is not a homomorphic image of $Sym(\mathbb{N})$.) We shall make use
of the following theorems of Serre.

\begin{Thm}{\cite[Theorem 16]{se}} \label{T:SL}
$SL(3,\mathbb{Z})$ has property (FA).
\end{Thm}

\begin{Thm}{\cite[Proposition 27]{se}} \label{T:nil}
Let $G = \langle g_{1}, \dots ,g_{n} \rangle$ be a finitely generated 
nilpotent group acting without inversion on the tree $T$. Suppose that
for each $1 \leq i \leq n$, there exists $t_{i} \in T$ such that
$g_{i}(t_{i}) = t_{i}$. Then there exists $t \in T$ such that
$g(t) = t$ for all $g \in G$.
\end{Thm}

For the rest of this section, let $Sym(\mathbb{N})$ act without inversion on
the tree $T$.

\begin{Lem} \label{L:sym1}
If $\pi \in Sym(\mathbb{N})$ contains no infinite cycles, then there exists
$t \in T$ such that $\pi(t) = t$.
\end{Lem}

\begin{proof}
There exists a sequence $\langle G_{n} \mid n \in \mathbb{N} \rangle$ of
nontrivial finite cyclic groups such that
$\pi \in \prod_{n} G_{n} \leqslant Sym(\mathbb{N})$. By Bass \cite{ba},
whenever the profinite group $\prod_{n} G_{n}$ acts without inversion
on a tree $T$, then for every $g \in \prod_{n} G_{n}$ there exists
$t \in T$ such that $g(t) = t$.
\end{proof}

\begin{Lem} \label{L:sym2}
Suppose that $\pi \in Sym(\mathbb{N})$ contains infinitely many infinite
cycles and no nontrivial finite cycles. Then there exists $t \in T$
such that $\pi(t) = t$.
\end{Lem}

\begin{proof}
Consider the left regular action of $SL(3, \mathbb{Z})$ on itself.
Let $h \in SL(3, \mathbb{Z})$ be an element of infinite order, and
let $H = \langle h \rangle$. Then $H$ has infinite index in
$SL(3, \mathbb{Z})$. Hence in the left regular action on
$SL(3, \mathbb{Z})$, $h$ contains infinitely many infinite cycles
and no finite cycles. Consequently if
$\Omega = \{ n \in \mathbb{N} \mid \pi(n) \ne n \}$, then there exists
a subgroup $G$ of $Sym(\mathbb{N})$ such that
\begin{enumerate}
\item $\pi \in G \leqslant Sym(\Omega) \leqslant Sym(\mathbb{N})$; and
\item the permutation group $\left( G, \Omega \right)$ is isomorphic
to the left regular action of $SL(3,\mathbb{Z})$ on itself.
\end{enumerate}
By Theorem \ref{T:SL}, there exists $t \in T$ such that $g(t) = t$
for all $g \in G$.
\end{proof}

\begin{Lem} \label{L:sym3}
Suppose that $\pi \in Sym(\mathbb{N})$ contains finitely many infinite
cycles and no nontrivial finite cycles, and that $\pi$ fixes 
infinitely many $n \in \mathbb{N}$. Then there exists $t \in T$ such
that $\pi(t) = t$.
\end{Lem}

\begin{proof}
There exist $\phi_{1}$, $\phi_{2} \in Sym(\mathbb{N})$ such that the
following conditions are satisfied.
\begin{enumerate}
\item $\phi_{1}$ and $\phi_{2}$ both contain infinitely many infinite
cycles and no nontrivial finite cycles.
\item $[\phi_{1}, \phi_{2}] = 1$.
\item $\pi = \phi_{1} \phi_{2}$.
\end{enumerate}
By Lemma \ref{L:sym2} and Theorem \ref{T:nil}, there exists $t \in T$
such that $g(t) = t$ for all $g \in G = \langle \phi_{1}, \phi_{2}
\rangle$.
\end{proof} 

\begin{proof}[Proof of Theorem \ref{T:amal}]
Let $\pi \in Sym(\mathbb{N})$ be any element. We shall show that $\pi$ fixes a
vertex of $T$. Express $\pi = \phi \psi$ as a product of disjoint
permutations such that $\psi$ has no infinite cycles and $\phi$ has
no nontrivial finite cycles. By Lemma \ref{L:sym1} and Theorem \ref{T:nil},
it is enough to show that $\phi$ fixes a vertex of $T$. Suppose that
$\phi \ne 1$. By Lemma \ref{L:sym2}, we can assume that $\phi$ contains
only finitely many infinite cycles. Let $\theta = \phi^{2}$. Then
$\theta$ contains $\ell$ infinite cycles for some $2 \leq \ell \in
\mathbb{N}$. Hence there exist $\tau_{1}$, $\tau_{2} \in Sym(\mathbb{N})$
such that the following conditions are satisfied.
\begin{enumerate}
\item $\tau_{1}$ and $\tau_{2}$ both contain finitely many infinite
cycles and no nontrivial finite cycles.
\item $\tau_{1}$ and $\tau_{2}$ are disjoint permutations.
\item $\theta = \tau_{1} \tau_{2}$.
\end{enumerate}
By Lemma \ref{L:sym3} and Theorem \ref{T:nil}, $\theta = \phi^{2}$
fixes a vertex of $T$. By 6.3.4 \cite{se}, $\phi$ also fixes a
vertex of $T$.
\end{proof}

\end{document}